\documentclass[12pt,letterpaper]{amsart}
\usepackage[utf8]{inputenc}
\usepackage{amsmath}
\usepackage{amssymb}
\usepackage{amsfonts}
\usepackage{amsthm}
\usepackage{tikz}
\usepackage{url}
\usepackage{hyperref}
\usepackage{cleveref}
\usepackage{subcaption}
\usepackage{stmaryrd}
\usepackage{wrapfig}
\usepackage[foot]{amsaddr}
\usepackage{enumitem}
\usepackage[normalem]{ulem}
\usepackage{mathtools}
\usepackage{fullpage}

\usepackage{ccicons}

\usetikzlibrary{matrix,decorations.pathreplacing}
\usetikzlibrary{positioning}

\newcommand{\bolda}{\mathbf{a}}

\newcommand{\bolde}{\mathbf{e}}

\newcommand{\boldx}{\mathbf{x}}
\newcommand{\boldy}{\mathbf{y}}
\newcommand{\boldalpha}{\boldsymbol{\alpha}}

\newcommand{\boldzero}{\mathbf{0}}
\newcommand{\boldone}{\mathbf{1}}
\newcommand{\inner}[2]{\langle \,#1\,,\, #2\, \rangle}
\newcommand{\biginner}[2]{\left\langle \,#1\,,\, #2\, \right\rangle}

\newcommand{\R}{\mathbb{R}}
\newcommand{\CC}{\mathbb{C}}

\newcommand{\adjp}{\nabla}
\newcommand{\facets}{\mathcal{F}}
\newcommand{\nodes}{\mathcal{V}}
\newcommand{\edges}{\mathcal{E}}

\newcommand{\digraph}[1]{\vec{#1}}
\newcommand{\red}[1]{\check{#1}}
\newcommand{\imag}{\mathbf{i}}

\newcommand{\rev}[1]{#1}

\DeclareMathOperator{\conv}{conv}
\DeclareMathOperator{\rank}{rank}
\DeclareMathOperator{\corank}{corank}

\newtheorem{theorem}{Theorem}
\newtheorem{proposition}[theorem]{Proposition}
\newtheorem{lemma}[theorem]{Lemma}
\newtheorem{corollary}[theorem]{Corollary}
\newtheorem{remark}[theorem]{Remark}

\theoremstyle{definition}
\newtheorem{definition}[theorem]{Definition}

\title{Facets and facet subgraphs of symmetric edge polytopes}
\author{Tianran Chen$^1$} 
\author{Robert Davis$^2$} 
\author{Evgeniia Korchevskaia$^{1,3}$}
\email{ti@nranchen.org}
\email{rdavis@colgate.edu}
\email{ekorchev@gatech.edu}
\address{$^1$Department of Mathematics, Auburn University Montgomery, Montgomery, AL, USA}
\address{$^2$Department of Mathematics, Colgate University, Hamilton, NY, USA}
\address{
    $^3$(Current affiliation) School of Mathematics, Georgia Institute of Technology, Atlanta, GA, USA
}
\thanks{
    TC and RD are supported by the National Science Foundation under grants no. 1923099
    and no. 1922998 respectively.
    TC and EK are supported by 
    Auburn University Montgomery (AUM) Grant-in-Aid Program.
    EK is also supported by the Undergraduate Research program 
    funded by AUM Department of Mathematics.
    This work is licensed under a CC-BY-NC-ND license
    \ccLogo~
    \ccAttribution~
    \ccNonCommercial~
    \ccNoDerivatives
    }

\begin{document}

\begin{abstract}
    Symmetric edge polytopes, a.k.a. PV-type adjacency polytopes,
    associated with undirected graphs have been defined and studied in several seemingly independent areas
    including number theory, discrete geometry, and dynamical systems. 
    In particular, the authors are motivated by the study of
    the algebraic Kuramoto equations of unmixed form whose Newton polytopes
    are the symmetric edge polytopes.
    
    The interplay between the geometric structure of symmetric edge polytopes
    and the topological structure of the underlying graphs
    has been a recurring theme in recent studies.
    In particular, ``facet/face subgraphs'' have emerged 
    as one of the central concepts in describing this symmetry.
    Continuing along this line of inquiry we provide a complete description
    of the correspondence between facets/faces of a symmetric edge polytope and
    maximal bipartite subgraphs of the underlying connected graph.
\end{abstract}



\maketitle

\section{Introduction}

For a simple connected graph $G$ with nodes $\nodes(G) = \{1,\dots,N\}$
and edge set $\edges(G)$,
its \emph{symmetric edge polytope} \cite{Matsui2011Roots}
(a.k.a. PV-type \emph{adjacency polytope}~\cite{Chen2019Unmixing})
is the convex polytope
$\conv \{ \pm(\bolde_i - \bolde_j) \mid \{i,j\} \in \edges(G) \}$ where 
$\bolde_0 = \boldzero$ and $\bolde_i$ is the $i^{th}$ standard basis.
In the context of Kuramoto models~\cite{Kuramoto1975Self},
the geometric structure of such polytopes has been
instrumental in understanding the root count for algebraic Kuramoto equations \cite{ChenDavisNODY, ChenDavisMehta2018Counting,Kuramoto1975Self}.
In the broader context, 
they have been studied  by number theorists, combinatorialists, and discrete geometers
motivated by several seemingly independent problems 
\cite{DelucchiHoessly, higashitani2019ARITHMETIC, HKMInterlacing,  Matsui2011Roots, ohsugi2014, OhsugiShibata2012Smooth, Rodriguez2002}.
These viewpoints are consolidated in recent work
by D'Alì, Delucchi, and Micha{\l}ek~\cite{DAliDelucchiMichalek2022Many}
which, among other contributions, sheds new light on the
structure of symmetric edge polytopes of bipartite graphs,
cycles, wheels, and graphs
consisting of two subgraphs sharing a single edge.
Using Gr\"obner basis methods, they provided explicit formulae for 
the number of facets and the volume 
of the symmetric edge polytopes associated with several classes of graphs. 

One recurring theme in these recent works
is the symmetry between the geometric structure of symmetric edge polytopes
and topological structure of the underlying graphs.
In particular, the concept of ``facet/face subgraphs'' is defined and studied
\cite{Chen2019Directed,DAliDelucchiMichalek2022Many}.
The present work is a continuation along this line of research.
The main contributions of this paper include the following descriptions
of the correspondence between faces of a symmetric edge polytope and
face subgraphs of the underlying connected simple graph:
\begin{itemize}
    \item
        Connected face subgraphs are exactly the 
        maximal bipartite subgraphs in their corresponding induced subgraphs.
        (\Cref{thm: facets are max bipartite} part (1))
    \item 
        Facet subgraphs are exactly the
        maximal bipartite subgraphs.
        (\Cref{thm: facets are max bipartite} part (2))
    \item
        The map from facets to facet subgraphs
        is surjective but not injective:
        each facet subgraph corresponds to an equivalence class
        of facets and each can be described as an assignment of edge orientations
        for the cut-set induced by the bipartite facet subgraph.
        From the view point of cut-set vectors,
        we provide a complete description of the equivalence class of facets
        corresponding to a given facet of the symmetric edge polytope,
        i.e., the fiber over a given facet subgraph.
        (\Cref{thm: facet parametrization})
    \item 
        We establish
        equivalences between geometric properties of faces and
        topological properties of the corresponding face subgraphs.
        (\Cref{thm: face properties})
    \item
        Armed with these results, we compute the number of facets
        of a symmetric edge polytope derived from a graph
        formed by joining an even and an odd cycle along a shared edge,
        generalizing a recent result of D'Alì, Delucchi, and Michałek
        \cite{DAliDelucchiMichalek2022Many}.
\end{itemize}

This paper is structured as follows.
\Cref{sec:prelim} states necessary definitions and notation.
\Cref{sec:adjpol} reviews the construction of symmetric edge polytopes
and face subgraphs.
\Cref{sec: facets review} provides a brief overview of
existing results on the interplay between facets/faces of
symmetric edge polytopes and their corresponding subgraphs.
Then, in \Cref{sec: face subgraphs,sec: facet constraints,sec: cut-set parametrization,sec: properties},
we develop the main results.
\Cref{sec:examples} illustrates how these results apply to 
a concrete class of non-bipartite graphs.
\Cref{sec: applications} highlights the important implications of the 
results in the study of Kuramoto equations
from dynamical systems and electrical engineering.

\section{Preliminaries and notation}\label{sec:prelim}

All graphs we consider will be simple.
For a graph $G$, let $\nodes(G)$ and $\edges(G)$ denote its sets of nodes and edges respectively.
We say $G$ is \emph{trivial} if $|\nodes(G)| \le 1$.
A graph $H$ is a \emph{subgraph} of $G$,
and use the notation $H \le G$,
if $\nodes(H) \subseteq \nodes(G)$ and $\edges(H) \subseteq \edges(G)$.
The edge connecting $i$ and $j$ is denoted
$\{i,j\}$ or $i \leftrightarrow j$.
A graph is \emph{bipartite} if it
is $2$-colorable; equivalently, if it has no odd cycles.
By a \emph{maximal} bipartite subgraph $H$ of $G$,
we mean a subgraph of $G$ that is bipartite and inclusion-maximal.
Such subgraphs are necessarily connected and \emph{spanning}
(i.e., $\nodes(H) = \nodes(G)$).
For a subset $V \subseteq \nodes(G)$,
the induced subgraph $G[V]$ is the subgraph
consisting of all edges $\{i,j\} \in \edges(G)$ where both $i,j \in V$.
With respect to a spanning tree $T$ of a connected graph $G$,
a \emph{fundamental cycle} of $G$ is the unique cycle formed by
an edge outside $T$ and a path in $T$.
For a digraph $\digraph{G}$,
the arrowhead emphasizes the distinction between $\digraph{G}$ and its underlying undirected graph $G$.
A directed edge, called an \emph{arc}, from $i$ to $j$ is denoted $(i,j)$,
and we use the notation $(i,j)^{-1} = (j,i)$.
The \emph{converse} of $\digraph{G}$, which reverses the orientation of all its edges, is denoted $\digraph{G}^{-1}$.


    A \emph{point configuration} $X = \{ \boldx_1,\dots,\boldx_m \} \subset \R^n$
    is a finite collection of distinct points.
    \rev{Its} (affine) \emph{dimension} is the dimension of the smallest affine subspace containing $X$.
    A nonempty \emph{face} of $X$ is a subset of $X$ for which 
    a linear functional $\inner{\cdot}{\boldalpha}$ is minimized.
    In this case, $\boldalpha$ is an \emph{inner normal} of the face. The emptyset
    $\varnothing$ is also a face of $X$. 
    Note that faces are themselves point configurations. The
    $0$-dimensional faces are \emph{vertices}
    and maximal proper faces are \emph{facets}.
    The set of all facets of $X$ is denoted by $\facets(X)$. 
    We say $X$ is \emph{ (affinely) dependent} if there are
    $\lambda_1,\dots,\lambda_m \in \R$, not all zero, with $\sum_{i=1}^m \lambda_i = 0$
    such that $\sum_{i=1}^m \lambda_i \boldx_i = \boldzero$.
    Otherwise, it is \emph{(affinely) independent}.
    Call $X$ a \emph{simplex} if $|X| = \dim(X) + 1$
    and a \emph{circuit} if it is an inclusion-minimal dependent set.
    Its \emph{corank}\footnote{%
        Here, the term ``corank'' is not to be confused with
        the different possible notions of corank in matroid theory --
        what we are calling ``corank'' is what a matroid theorist would call ``nullity''.
    }
    is $|X| - \dim(X) - 1$; 
    Its convex hull
    $\conv(X) = 
    \left\{
        \sum_{i=1}^m \lambda_i\boldx_i \mid \lambda_1,\dots,\lambda_m \geq 0, \, \sum_{i=1}^m \lambda_i = 1\right
    \}
    $ is a \emph{convex polytope}.
\section{Symmetric edge polytopes}\label{sec:adjpol}

For a connected graph $G$ with nodes $\nodes(G) = \{1,\dots,N\}$,
its \emph{symmetric edge polytope~\cite{DAliDelucchiMichalek2022Many,Matsui2011Roots}}
is the convex polytope 
$\conv \{ \pm(\bolde_{i} - \bolde_{j}) \mid \{i,j\} \in \edges(G) \} \subset \R^N$,
where
$\bolde_i \in \R^N$ is the vector with 1 in the $i$-th entry and zero elsewhere.
Since we are mostly interested in combinatorial aspects of symmetric edge polytopes,
it is more convenient to focus on the underlying point configuration.
We define
\begin{align*}
    \red{\adjp}_G &=
    \{ \pm( \bolde_{i-1} - \bolde_{j-1}) \mid \{i,j\} \in \edges(G) \}
    \subset \R^n = \R^{N-1}
    \quad\text{and}
    \\
    \bar{\adjp}_G &=
    \{ \pm( \bolde_i - \bolde_j) \mid \{i,j\} \in \edges(G) \}
    \subset \R^{N}
\end{align*}
with the convention that $\bolde_0 = \boldzero$.
These two point configurations have the same intrinsic geometric properties,
and they only differ in the ambient space in which they are embedded:
$\red{\adjp}_G$ is a full-dimensional point configuration in $\R^{N-1}$
whereas $\bar{\adjp}_G$ is a codimension-1 point configuration in $\R^N$,
and $\red{\adjp}_G$ is precisely the projection of $\bar{\adjp}_G$
onto the last $N-1$ coordinates.
The check mark notation in $\red{\adjp}_G$ is a reminder
that it is a projection to a lower-dimensional subspace.
We also extend this notation to their subsets,
e.g., we identify any subset $X \subseteq \bar{\adjp}_G$
with its projection $\red{X} \subseteq \red{\adjp}_G$.
When describing subsets of $\red{\adjp}_G$ or $\bar{\adjp}_G$,
the ``codimension'' of a subset always refers to the codimension
relative to $\red{\adjp}_G$ or $\bar{\adjp}_G$ themselves,
regardless of the ambient space.
The projection that maps $\bar{\adjp}_G$ to $\red{\adjp}_G$ 
is a unimodular equivalence between the two configurations,
\rev{%
    and the two are simply different embeddings of the same polytope
    into the Euclidean space.
}%
The distinction between the two will not be relevant when referencing intrinsic geometric properties, in which case we will simply use $\adjp_G$.


\rev{
In Ref.~\cite{Higashitani2015Smooth2},
Higashitani extended this construction to digraphs}:
for a digraph $\digraph{G}$, we define
\begin{align*}
    \red{\adjp}_{\digraph{G}} &= 
    \{ \bolde_{i-1} - \bolde_{j-1} \mid (i,j) \in \edges(\digraph{G}) \} \subset \R^{N-1}
    \quad\text{and}
    \\
    \bar{\adjp}_{\digraph{G}} &= 
    \{ \bolde_{i} - \bolde_{j} \mid (i,j) \in \edges(\digraph{G}) \} \subset \R^N.
\end{align*}
Here, $\bolde_i - \bolde_j \in \bar{\adjp}_{\digraph{G}}$ no longer implies 
$\bolde_j - \bolde_i \in \bar{\adjp}_{\digraph{G}}$, 
thus $\conv(\bar{\adjp}_{\digraph{G}})$ may not be a symmetric edge polytope.
The notation 
$\red{\adjp}_G,\bar{\adjp}_G$ and 
$\red{\adjp}_{\digraph{G}},\bar{\adjp}_{\digraph{G}}$
extend naturally to subgraphs of $G$ and 
subdigraphs of a digraph $\digraph{G}$, respectively, by restriction.

By construction, $\boldzero$ is always an interior point of $\conv(\red{\adjp}_G)$,
which allows the inner normals to be normalized to the following certain form.

\begin{lemma}\label{lem:facet-def}
    For a connected nontrivial graph $G$, 
    a nonempty subset $F \subsetneq \red{\adjp}_G$
    is a face if and only if there is
    a nonzero vector $\red{\boldalpha} \in \R^n$
    such that
    \begin{align*}
        \inner{ \boldx }{ \red{\boldalpha}} &=   -1
        \quad\text{for any } \boldx \in F,\; \text{and}\\
        \inner{ \boldx   }{ \red{\boldalpha}} &> -1
        \quad\text{for any } \boldx \in \red{\adjp}_G \setminus F
    \end{align*}
\end{lemma}

The same description extends to inner normals for faces of $\bar{\adjp}_G$ since 
each inner normal of $\bar{\adjp}_G$ projects down to 
an inner normal of $\red{\adjp}_G$, and
each inner normal of $\red{\adjp}_G$ lifts to 
an equivalence class of inner normals of $\bar{\adjp}_G$.

\section{Facets, faces, and associated subgraphs}\label{sec: facets}


Faces and facets of a symmetric edge polytope
have been studied from several different viewpoints \cite{DAliDelucchiMichalek2022Many}.
We continue this line of inquiry through a graph-theoretical approach.
The central theme of this paper is the interplay between 
combinatorial properties of faces of $\adjp_G$ and 
graph-theoretic properties of subgraphs of $G$
through the connection of face subgraphs.
Throughout this section,
we fix $G$ to be a nontrivial, connected, and simple graph.
    
\begin{definition} \label{def:subgraphs}
    For a nonempty subset $X \subseteq \bar{\adjp}_G$,
    we define $\digraph{G}_X$ and $G_X$ to be the subgraphs with node and edge sets 
    \begin{align*}
        \nodes(\digraph{G}_X) &= \{ 
            i \mid 
            \bolde_i - \bolde_j \in X
            \text{ or }
            \bolde_j - \bolde_i \in X
            \text{ for some } j
        \}
        \\
        \edges(\digraph{G}_X) &= \{ \; (i,j) \;\mid\; \bolde_i - \bolde_j \in X \}
        \\
        \nodes(G_X) &= \{ 
            i \mid 
            \bolde_i - \bolde_j \in X
            \text{ or }
            \bolde_j - \bolde_i \in X
            \text{ for some } j
        \}
        \\
        \edges(G_X)           &= \{ \; 
            \{ i, j \} 
            \;\mid\; \bolde_i - \bolde_j \in X \text{ or } \bolde_j - \bolde_i \in X 
        \},
    \end{align*}
    respectively.
    If $F$ is a face (resp. facet) of $\adjp_G$, 
    then $\digraph{G}_F$ is the \emph{face (resp. facet) subdigraph} 
    associated with $F$,
    and $G_F$ is the associated \emph{face (resp. facet) subgraph}.
\end{definition}

These concepts are defined and studied in recent works
\cite{Chen2019Directed,DAliDelucchiMichalek2022Many}.
The notational conventions are chosen to mirror the connection
between $\digraph{G}$ and $\adjp_{\digraph{G}}$,
and they interact in an expected way:
for any $X \subseteq \adjp_G$, we have $\adjp_{\digraph{G}_X} = X$, and 
for any subgraph $\digraph{H} \le \digraph{G}$, 
we have $\digraph{G}_{\adjp_{\digraph{H}}} = \digraph{H}$.

Note that points in $X$
are exactly the columns in the incidence matrix of $\digraph{G}_X$,
which will be denoted by $Q(\digraph{G}_X)$.
The \emph{truncated incidence} matrix $\red{Q}(\digraph{G}_X)$,
obtained by removing the first row of $Q(\digraph{G}_X)$,
corresponds to points in the projection $\red{X} \subset \red{\adjp}_G$.


Face and facet subgraphs 
were studied in Refs.~\cite{DAliDelucchiMichalek2022Many,GordonPetrov2017,higashitani2019ARITHMETIC}.
\Cref{sec: facets review} reviews recent results
on the interplay between facets and their corresponding subgraphs.
Then, in \Cref{sec: face subgraphs,sec: facet constraints,sec: cut-set parametrization,sec: properties},
we develop the main results.
In particular, we provide a complete description of facet and face subgraphs
as well as the equivalence classes of facet subdigraphs.
Implications of these results in algebraic Kuramoto equations
are highlighted in \Cref{sec: applications}.

\subsection{Recent results on facets and facet subgraphs}\label{sec: facets review}

Facets of $\adjp_G$ associated with even cycles
are described from the viewpoint of Ehrhart theory
by Ohsugi and Shibata~\cite{OhsugiShibata2012Smooth}
and from the viewpoint of Lipschitz polytope
by Gordon and Petrov \cite{GordonPetrov2017}.
Explicit descriptions of the facets of $\adjp_G$,
when $G$ is a tree or cycle, are also established~\cite{ChenDavisMehta2018Counting}.
Using Gr\"obner basis methods,
D'Alì, Delucchi, and Micha{\l}ek \cite{DAliDelucchiMichalek2022Many} provided
more detailed descriptions for the faces of $\adjp_G$.
In particular, they showed that
unimodular simplices contained in a facet $F \in \facets(\adjp_G)$
correspond to spanning trees of $G_F$
\rev{\cite[Corollary 3.3]{higashitani2019ARITHMETIC}},
and for a connected bipartite graph $G$,
the total number of facets is bounded by $2^{|\nodes(G)|-1}$.
This bound is exact if $G$ is a tree.
Indeed, in this special case, the numbers of faces of any dimension
are offered by Delucchi and Hoessly \cite{DelucchiHoessly}.
In \Cref{sec: cut-set parametrization}, we provide graph-theoretic refinements
to these results.

For an even cycle $C_{2k}$, the number of faces
of each dimension of $\adjp_{C_{2k}}$,
is also computed \rev{in} \cite[\rev{Proposition 30}]{DAliDelucchiMichalek2022Many}.
Moreover, if $G_1$ and $G_2$ are both connected bipartite graphs,
and $G$ is formed by identifying an edge of $G_1$ with an edge of $G_2$, then
it has been shown that $| \facets(\adjp_G) | = \frac{1}{2} f_1 f_2$
where $f_1 = | \facets(\adjp_{G_1}) |$ 
and   $f_2 = | \facets(\adjp_{G_2}) |$ 
\cite[\rev{Proposition 37}]{DAliDelucchiMichalek2022Many}.
This result can be applied recursively and extended to 
graphs formed by joining multiple even cycles consecutively by an edge
\cite[\rev{Corollary 38}]{DAliDelucchiMichalek2022Many}.
We will extend these to non-bipartite graphs.

\subsection{Characterizing face and facet subgraphs}\label{sec: face subgraphs}

Our goal is to clarify the structure of the map $F \mapsto G_F$ 
between facets and faces of $\adjp_G$ and subgraphs of $G$.
    We first show that connected face and facet subgraphs 
    associated with faces and facets of $\adjp_G$ are exactly the maximal bipartite subgraphs of 
    induced subgraphs of $G$ (\Cref{thm: facets are max bipartite}).
    In particular, the map $F \mapsto G_F$ is a surjective map
    from $\facets(\adjp_G)$ to the set of maximal bipartite subgraphs of $G$.
    \Cref{cor: face subgraphs components} generalizes this description
    to components of face subgraphs.

\begin{theorem}\label{thm: facets are max bipartite}
    Let $H$ be a nontrivial connected subgraph of $G$.
    \begin{enumerate}
        \item
            $H$ is a face subgraph of $G$ 
            if and only if it is a maximal bipartite subgraph of $G[\nodes(H)]$.
        \item
            $H$ is a facet subgraph of $G$ if and only if it is a maximal bipartite subgraph of $G$.
    \end{enumerate}
\end{theorem}

    Note that part (2) of this theorem can be derived from
    \cite[Theorem 3.1]{HKMInterlacing},
    which provides a general description of facet-defining labels
    for nodes in $G$.
    Here, we provide an alternative derivation
    as a special case of part (1).
\begin{proof}
    First, note that facet subgraphs and maximal bipartite subgraphs 
    are necessarily connected, nontrivial, and spanning,
    therefore part (2) is a special case of part (1).
    If $F \ne \varnothing$ is a proper face of $\bar{\adjp}_G$, 
    then, by \Cref{lem:facet-def}, there exists an 
    $\boldalpha \in \{ \boldone \}^\perp \subset \R^N$  such that
    \begin{align*}
        \biginner{ \boldx }{ \boldalpha } &= -1
        \quad\text{for all } \boldx \in F
        &&\text{and} &
        \biginner{ \boldx }{ \boldalpha } &> -1
        \quad\text{for all } \boldx \in \bar{\adjp}_G \setminus F.
    \end{align*}
    To show $G_F$ is bipartite, suppose
    $i_1 \leftrightarrow \cdots \leftrightarrow i_\ell \leftrightarrow i_1$
    is a cycle in $G_F$. 
    Then there are $\lambda_1,\dots,\lambda_\ell \in \{ \pm 1 \}$ such that
    $\lambda_j (\bolde_{i_j} - \bolde_{i_{j+1}}) \in F$ for $j=1,\dots,\ell$
    with $i_{\ell+1} = i_1$.
    By the above,
    \begin{align*}
        \biginner{ \lambda_j(\bolde_{i_j} - \bolde_{i_{j+1}}) }{ \boldalpha } &= -1 &
        &\text{i.e.,} &
        \biginner{ \bolde_{i_j} - \bolde_{i_{j+1}} }{ \boldalpha } &= 
        -\lambda_j &
        &\text{for } j = 1,\dots,\ell.
    \end{align*}
    Summing both sides over all $j$ produces
    \[
        0 =
        \biginner{ \boldzero }{ \boldalpha } =
        \biginner{ \sum_{j=1}^\ell \bolde_{i_j} - \bolde_{i_{j+1}} }{ \boldalpha } = 
        - \sum_{j=1}^\ell \lambda_j.
    \]
    Since $\lambda_j \in \{ \pm 1 \}$,
    $\ell$ must be even,
    i.e., any cycle in $G_F$ must be even,
    and $G_F$ is bipartite.
    
    To show $G_F$ is a maximal bipartite subgraph of $G[\nodes(G_F)]$,
    consider a spanning tree $T$ of $G_F$,
    which is necessarily a spanning tree of $G[\nodes(G_F)]$.
    If $B$ is a bipartite graph such that $G_F < B \le G[\nodes(G_F)]$,
    then any edge $\{i,i'\} \in \edges(B) \setminus \edges(G_F)$ is also outside $T$.
    The fundamental cycle formed by $\{i,i'\}$ and the unique path
    $i = i_1 \leftrightarrow \cdots \leftrightarrow i_\ell = i'$ in $T$
    is contained in $B$ and hence must be an even cycle.
    That is, $\ell$ is even.
    Since the path $i_1 \leftrightarrow \cdots \leftrightarrow i_\ell$ 
    is in $T \le G_F$,
    there are $\lambda_1,\dots,\lambda_{\ell-1} \in \{ \pm 1 \}$ such that
    $\lambda_j (\bolde_{i_j} - \bolde_{i_{j+1}}) \in F$ for $j=1,\dots,\ell-1$.
    As in the paragraph above,
    \[
        \biginner{ \bolde_{i} - \bolde_{i'} }{ \boldalpha } = 
        \biginner{ \sum_{j=1}^{\ell-1} (\bolde_{i_j} - \bolde_{i_{j+1}}) }{ \boldalpha } =
        - \sum_{j=1}^{\ell-1} \lambda_j.
    \]
    Also, by \Cref{lem:facet-def},
    $
        \biginner{ \pm(\bolde_{i} - \bolde_{i'}) }{ \boldalpha } = 
        \mp \sum_{j=1}^{\ell-1} \lambda_j 
        \ge -1.
    $
    Recall that $\lambda_{j} \in \{ \pm 1 \}$ and $\ell$ is even.
    So
    either $\biginner{ +(\bolde_{i} - \bolde_{i'}) }{ \boldalpha }$
    or     $\biginner{ -(\bolde_{i} - \bolde_{i'}) }{ \boldalpha }$
    must be $-1$, and hence
    either $+(\bolde_{i} - \bolde_{i'})$ 
    or     $-(\bolde_{i} - \bolde_{i'})$
    is also contained in $F$.
    This implies $\{i,i'\} \in \edges(G_F)$,
    contradicting our assumption.
    Therefore, $G_F$ is a not contained in any larger
    bipartite subgraphs of $G[\nodes(G_F)]$.
    
    For the converse, suppose a connected subgraph $B \le G$ 
    is a maximal bipartite subgraph of $G[\nodes(B)]$
    with the partition $\nodes(B) = V_+ \cup V_-$.
    Define $\boldalpha = (\alpha_1,\dots,\alpha_N)$ with
    \begin{align*}
        \alpha_i &=
        \begin{cases}
            +1/2 &\text{if } i \in V_+ \\
            -1/2 &\text{if } i \in V_- \\
            0    &\text{otherwise}.
        \end{cases}
    \end{align*}
    Then 
    \[
        \inner{ \bolde_i - \bolde_j }{ \boldalpha } =
        \begin{cases}
            0     &\text{if } i,j \in V_+ \text{ or } i,j \in V_- \text{ or } i,j \not\in V_+ \cup V_- \\
            \pm 1 &\text{if } i \in V_\pm \text{ and } j \in V_\mp \\
            \pm 1/2 &\text{if exactly one of $i,j$ is in $V_+ \cup V_-$}        \end{cases}
    \]
    for any $i,j \in \nodes(G)$.
    In particular,  
    since $V_+$ and $V_-$ partition the nodes of $B$, which is bipartite,
    $\inner{ \bolde_i - \bolde_j }{ \boldalpha } = \pm 1$
    for every $\{i,j\} \in \edges(B)$.
   
    Let
    \[
        F = \{ 
            \bolde_i - \bolde_j \mid 
            \inner{ \bolde_i - \bolde_j }{ \boldalpha } = -1
        \}.
    \]
    The linear functional $\inner{ \,\cdot\, }{ \boldalpha }$
    attains its minimum over $\adjp_G$ on $F$,
    and thus $F$ is a face of $\adjp_G$.
    Moreover, $G_F$ is a bipartite subgraph of $G[\nodes(B)]$ that contains $B$.
    But $B$ is assumed to be a maximal bipartite subgraph of $G[\nodes(B)]$,
    so $B = G_F$, i.e., $B$ is a face subgraph.
\end{proof}

\begin{wrapfigure}[10]{r}{0.25\textwidth}
    \centering
    \begin{tikzpicture}[scale=1.25,
        every node/.style={circle,thick,draw,inner sep=2pt}]
        \node (1) at (   0,  1) {1};
        \node (2) at (   1,  0) {2};
        \node (3) at (   0, -1) {3};
        \node (4) at (  -1,  0) {4};
        \path[thick,-latex] (1) edge (2);
        \path[thick,-latex] (4) edge (3);
        \path[gray,densely dotted] (1) edge (4);
        \path[gray,densely dotted] (2) edge (3);
    \end{tikzpicture}
    \caption{
        A disconnected face subgraph.
    }
    \label{fig: nonmaximal face subgraph}
\end{wrapfigure}
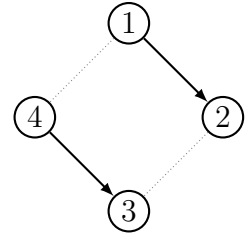
The connectedness condition in \Cref{thm: facets are max bipartite}
is important as disconnected face subgraphs may not be 
maximal bipartite subgraph of their associated induced subgraphs.
For example, \Cref{fig: nonmaximal face subgraph} shows a face subgraph
associated with the 4-cycle $G = C_4$,
which is not a maximal bipartite subgraph of its associated induced subgraph
$G[\{1,2,3,4\}] = G$.
Nonetheless, the argument still applies
to individual connected components of a face subgraph.
From this observation we can derive the following generalization.

\begin{corollary}\label{cor: face subgraphs components}
    Let $F$ be a face of $\adjp_G$.
    Then each connected component $H$ of $G_F$ is a maximal bipartite subgraph of $G[\nodes(H)]$.
    \qed
\end{corollary}

The proof of \Cref{thm: facets are max bipartite} also shows that
for a maximal bipartite subgraph,
there is a canonical choice of edge orientations
that defines a facet subdigraph and hence a facet.

\begin{corollary}\label{cor: canonical facet}
    For a maximal bipartite subgraph $B \le G$
    with the partition $\nodes(B) = V_+ \cup V_-$,
    \[
        F =
        \{ \bolde_i - \bolde_j \mid \{i,j\} \in \edges(B), i \in V_- \text{ and } j \in V_+ \}
    \]
    is a facet of $\bar{\adjp}_G$, defined by the facet inner normal
    $\boldalpha = (\alpha_1,\dots,\alpha_N)^\top$ with
    \[
        \alpha_i = 
        \begin{cases}
            + 1/2 &\text{if } i \in V_+ \\
            - 1/2 &\text{if } i \in V_- \, ,
        \end{cases}
    \]
    and the associated facet subdigraph is the digraph with edges set
    \[
        \edges(\digraph{G}_F) =
        \{ (i,j) \mid \{i,j\} \in \edges(B), \, i \in V_-, \, j \in V_+ \}.
    \]
\end{corollary}

The choice of the labels $V_+$ and $V_-$ is, of course, arbitrary:
permuting the two will result in the facet $-F$ defined by $-\boldalpha$
associated with the facet subdigraph $\digraph{B}^{-1}$.
This partition defines a 
unique \emph{cut-set},
which is the set of edges that go across the partition.
The construction in the above corollary
can be interpreted as a special type of cut-set.

\begin{remark}\label{rmk: canonical orientation}
    For a maximal bipartite subgraph $B \le G$,
    the canonical choice of edge orientations that produces $\digraph{B}$
    in \Cref{cor: canonical facet} are exactly the edge orientation assignments that
    ensure the cut-set (which includes all edges in $B$)
    has a uniform direction across the cut $(V_+,V_-)$ of $B$,
    i.e., all arcs of $\digraph{B}$ are from $V_-$ to $V_+$.
    Such an assignment will be referred to as a
    \emph{canonical edge orientation} for $B$
    as well as its spanning subgraphs.
\end{remark}

\subsection{Cyclic constraints on facet subdigraphs}\label{sec: facet constraints}

\Cref{thm: facets are max bipartite} shows that the map $F \mapsto G_F$
is a surjective map from $\facets(\adjp_G)$ to the set of maximal bipartite subgraphs of $G$.
However, this map is not injective.
Indeed, as clarified in \Cref{cor: canonical facet},
there is at least a pair of canonical choices of facets $F$ and $-F$ 
associated with a given maximal bipartite subgraph $B$ of $G$ which correspond to cut-sets having either of the two uniform directions.
In general, if we let $V_+ \cup V_- = \nodes(B) = \nodes(G)$ be the partition in $B$,
then each facet $F \in \facets(\adjp_G)$ such that $G_F = B$
can be described as an assignment of edge orientations
for the cut-set induced by the cut $(V_+,V_-)$ in $B$.
However, not every such assignment will produce a facet of $\adjp_G$.
The constraints on such assignments has been studied in
\cite[Section 6]{GordonPetrov2017} from the viewpoint of metric spaces.
In the following, we describe constraints
on the possible choices that will result in facets 
in terms of \emph{oriented cycles}: directed
cycles in which each node is the head of
exactly one arc and the tail of exactly one arc.

\begin{theorem} \label{thm: facets have half and half}
        If $F$ is a facet of $\adjp_G$,
        then for any cycle $\digraph{O}$ in $G$ with an assigned orientation,
        \[ 
            |\edges(\digraph{G}_F) \cap \edges( \digraph{O})| = 
            |\edges(\digraph{G}_F) \cap \edges(\digraph{O}^{-1})|.
        \]
        
\end{theorem}

Note that $G_F$ being a maximal bipartite subgraph
already implies that $|\edges(G_F) \cap \edges(O)|$ is even
for any cycle $O$ in $G$.
This theorem states that $\edges(\digraph{G}_F)$
consists of two subsets of arcs of equal size having opposite orientations.
Later in this section, we will show 
that, under an additional dimensional condition, the converse is also true.

\begin{proof}
    
    
            
            
        Recall that the incidence matrix $Q(\digraph{G}_F)$ is totally unimodular
        \cite[Lemma 2.6]{Bapat2010Graphs},
        and the reduced inner normal $\check{\boldalpha}$, being the vector satisfying 
        $\check{\boldalpha}^\top \check{Q}(\digraph{G}_F) = -\boldone^\top$,
        must be an integer vector.
        Then $\boldalpha = (0,\check{\boldalpha})$ is also an integer vector,
        and $\pm \inner{\bolde_i-\bolde_j}{\boldalpha} > -1$ is an integer
        for any $\bolde_i - \bolde_j \not\in F$.
        This implies that $\inner{\bolde_i-\bolde_j}{\boldalpha}=0$
        for any $\bolde_i - \bolde_j \not\in \pm F$.
        
        Suppose $G$ contains a cycle $O$ of length $m$ having edges 
        $i_1 \leftrightarrow i_2 \leftrightarrow \cdots \leftrightarrow i_{m+1}$
        with $i_{m+1} = i_1$,
        and let $E = \edges(G_F) \cap \edges(O)$.
        Then $\boldzero = \sum_{r=1}^{m} \rev{(\bolde_{i_r}-\bolde_{i_{r+1}})}$ implies
        \begin{align*}
            0 &=
            \biginner{\sum_{r=1}^{m}\rev{(\bolde_{i_r}-\bolde_{i_{r+1}}) } }{\boldalpha} \\  &=
            \sum_{\{j,j+1\} \in E } \biginner{ \bolde_{i_{r_j}}-\bolde_{i_{r_j+1}} }{\boldalpha} +
            \sum_{\{j,j+1\} \in \edges(O) \setminus E } \biginner{ \bolde_{i_{r_j}}-\bolde_{i_{r_j+1}} }{\boldalpha} \\ &=
            \sum_{\{j,j+1\} \in E } \biginner{ \bolde_{i_{r_j}}-\bolde_{i_{r_j+1}} }{\boldalpha},
        \end{align*}    
        since $\inner{\bolde_i-\bolde_j}{\boldalpha}=0$
        for any $\bolde_i - \bolde_j \not\in \pm F$.
        Moreover, $\inner{\bolde_i-\bolde_j}{\boldalpha} = \pm 1$
        for any $\bolde_i - \bolde_j \in F$,
        i.e., each term in the above sum is $\pm 1$.
        We conclude that $|E| = | \edges(G_F) \cap \edges(O) |$ is even.
        
        Let $\digraph{O}$ with arcs  
        $i_1 \to i_2 \to \cdots \to i_{m+1}$
        be the corresponding oriented cycle.
        Then $\inner{ \bolde_{i_{r_j}}-\bolde_{i_{r_j}+1} }{\boldalpha} = -1$
        implies $(i_{r_j},i_{r_{j+1}}) \in \edges(\digraph{G}) \cap \edges(\digraph{O})$
        and $\inner{ \bolde_{i_{r_j}}-\bolde_{i_{r_j}+1} }{\boldalpha} = 1$
        implies $(i_{r_{j+1}},i_{r_j}) \in \edges(\digraph{G}) \cap \edges(\digraph{O}^{-1})$.
        Therefore,
        \[ 
            |\edges(\digraph{G}_F) \cap \edges( \digraph{O})| = 
            |\edges(\digraph{G}_F) \cap \edges(\digraph{O}^{-1})|.
            \qedhere
        \]
\end{proof}

\Cref{thm: facets have half and half} shows necessary conditions
for a subset $F$ of $\bar{\adjp}_G$ to be a facet:
the intersection between $\digraph{G}_F$ and any cycle of $G$ itself
must satisfy a balancing condition, i.e.,
$|\edges(\digraph{G}_F) \cap \edges( \digraph{O})| = 
 |\edges(\digraph{G}_F) \cap \edges(\digraph{O}^{-1})| $
for any oriented cycle $\digraph{O}$.
By itself, however, this condition is not sufficient to define a facet.
In the following, we will show that under additional assumptions,
this balancing condition is also a sufficient condition.

\begin{theorem}\label{thm: half and half make a facet}
    Let $F$ be a codimension 1 subset of $\adjp_G$
    such that $\boldzero \not\in \conv(F)$.
    Then $F$ is a facet of $\adjp_G$
    if, for any oriented cycle $\digraph{O}$ in $\digraph{G}$, 
    \[ 
        |\edges(\digraph{G}_F) \cap \edges( \digraph{O})| = 
        |\edges(\digraph{G}_F) \cap \edges(\digraph{O}^{-1})|\rev{.}
    \]
\end{theorem}

\begin{proof}
    It is sufficient to consider the embedding $\red{\adjp}_G \subset \R^n$
    and assume $F \subset \red{\adjp}_G$.
    Since $F$ is a codimension 1 subset of $\red{\adjp}_G$ that does not contain $\boldzero$
    in its convex hull,
    \[
        N - 2 = 
        \dim(F) = \rank(Q(\digraph{G}_F)) - 1 = |\nodes(\digraph{G}_F)| - k - 1
    \]
    where $Q(\digraph{G}_F)$ is the incidence matrix whose columns are points in $F$
    and $k$ is the number of weakly connected components in $\digraph{G}_F$
    \cite[Theorem 2.3]{Bapat2010Graphs}.
    Therefore $\digraph{G}_F$ is necessarily (weakly) connected and spanning,
    and $G_F$ is connected and spanning.
    
    If $T$ is a spanning tree of $G_F$,
    then $T$ is also a spanning tree of $G$.
    Let $\digraph{T}$ be the corresponding directed subgraph of $\digraph{G}_F$.
    By \rev{\cite[Corollary 3.2]{higashitani2019ARITHMETIC}}
    (\cite[\rev{Corollary 13}]{DAliDelucchiMichalek2022Many}),
    points in $\Delta = \red{\adjp}_{\digraph{T}} \subset \red{F}$ form a simplex,
    and $\check{Q}(\digraph{T})$ is nonsingular.
    Additionally, if $\check{\boldalpha}$ is the unique solution to
    $
        \red{\boldalpha}^\top \check{Q}(\digraph{T}) = - \boldone^\top,
    $
    then 
    \[
        \inner{ \bolde_i - \bolde_j }{ \red{\boldalpha} } = -1
        \quad
        \text{for all } (i,j) \in \digraph{T}.
    \]
    Since $F$ is assumed to be a codimension 1 subset in $\red{\adjp}_G$,
    i.e., $\dim(F) = \dim(\Delta)$,
    $F$ must be contained in the affine span of $\Delta$.
    Consequently,
    \[
        \inner{ \bolde_i - \bolde_j }{ \red{\boldalpha} } = -1
        \quad
        \text{for all } \bolde_i - \bolde_j \in F.
    \]
    Recall that $\boldzero$ is not contained in the convex hull of $F$,
    i.e., $F$ cannot contain both $\pm(\bolde_i - \bolde_j)$ for any $\{i,j\}$,
    so $F$ and $-F$ are disjoint.
    For any $\bolde_i - \bolde_j \in -F$ and hence outside $F$,
    it is clear that
    $\inner{ \bolde_i - \bolde_j }{ \red{\boldalpha} } = +1$.
    
    For any $\bolde_i - \bolde_j \in \red{\adjp}_G \setminus (F \cup (-F))$,
    the corresponding undirected edge $\{i,j\}$ is outside $G_F$ 
    and hence outside $T$.
    Consider the fundamental cycle $O$ formed by $\{i,j\}$ and the path 
    $i = i_1 \leftrightarrow \cdots \leftrightarrow i_m \leftrightarrow i_{m+1} = j$ in $T$.
    We have
    \[
        \bolde_i - \bolde_j = 
        \sum_{j=1}^m \rev{ (\bolde_{i_j} - \bolde_{i_{j+1}}) },
    \]
    and there are $\lambda_1,\dots,\lambda_m \in \{ \pm 1 \}$
    such that $\lambda_j (\bolde_{i_j} - \bolde_{\rev{ i_{j+1} } }) \in \Delta \subseteq F$.
    By the assumption that
    $
        |\edges(\digraph{G}_F) \cap \edges( \digraph{O})| = 
        |\edges(\digraph{G}_F) \cap \edges(\digraph{O}^{-1})|,
    $
    $m$ must be even, and $\lambda_1 + \cdots + \lambda_m = 0$.
    Therefore,
    \begin{align*}
        \inner{ \bolde_i - \bolde_j }{ \boldalpha } &=
        \sum_{j=1}^m \inner{ \bolde_{i_j} - \bolde_{i_{j+1}} }{ \boldalpha } = 
        \sum_{j=1}^m \lambda_j \inner{ \lambda_j(\bolde_{i_j} - \bolde_{i_{j+1}}) }{ \boldalpha }
        = 
        \sum_{j=1}^m \lambda_j (-1) = 0.
    \end{align*}
    That is, the linear functional $\inner{ \cdot }{\boldalpha}$
    takes the value of $-1$ on $F$, and it is nonnegative on
    $\adjp_G \setminus F$.
    Therefore $F$ is a facet.
\end{proof}

\subsection{Parameterizing facets with cut-set vectors}\label{sec: cut-set parametrization}

Facets of $\adjp_G$ correspond to maximal bipartite subgraphs of $G$ 
through the map $F \mapsto G_F$.
In general, multiple facets will be mapped to the same facet subgraph.
\Cref{thm: facets have half and half,thm: half and half make a facet}
gave necessary and sufficient conditions to identify facets in the fiber over a given facet subgraph in terms of oriented cycles. 
    In this section, we refine these constraints
    into ``independent'' equations
    and thereby
provide a complete description
of the equivalence class of facets corresponding to the same facet subgraph.
It will form the foundation for counting and generating facets of $\adjp_G$.

The description makes use of the fundamental cycle vectors and cut-set vectors.
For a facet subgraph $G_F$ and a spanning tree $T$ of $G_F$,
let $\digraph{T}$ be the corresponding subdigraph of $\digraph{G}_F$.
Since a facet subgraph is necessarily connected and spanning, $T$ is also a spanning tree of $G$.
Any edge $e \in \edges(G) \setminus \edges(T)$ induces a fundamental cycle $O$ with respect to $T$.
With an arbitrary choice of the orientation,
the oriented cycle $\digraph{O}$ can be expressed as an incidence vector
$\mathbf{c}_{\digraph{T}}(e) = (c_1,\dots,c_n)^\top$ with respect to 
the ordered list of arcs in $\digraph{T}$ so that $\digraph{e}^{-1}$ 
corresponds to the point $Q(\digraph{T}) \mathbf{c}_{\digraph{T}}(e) \in \red{\adjp}_G$.
    In other words,
    $\mathbf{c}_{\digraph{T}}(e) \in \{ +1, 0, -1 \}^{|\edges(\digraph{T})|}$,
    with each entry indicating whether the orientation of the corresponding edge
    of $\digraph{T}$ agrees or disagrees (or is not involved)
    with the orientation of $\digraph{O}$.
    
\begin{figure}[ht]
    \centering
    \begin{subfigure}{0.3\textwidth}
        \centering
        \begin{tikzpicture}[scale=1.0,
            every node/.style={circle,thick,draw,inner sep=2pt},
            every edge/.style={draw,thick}]
            \node (1) at (   0,  1) {1};
            \node (2) at (   1,  1) {2};
            \node (3) at (   1,  0) {3};
            \node (4) at (   0,  0) {4};
            \node (5) at (  -1,  0) {5};
            \node (6) at (  -1,  1) {6};
            \node (7) at (-0.5,1.75) {7};
            \path (1) edge (2);
            \path (2) edge (3);
            \path (3) edge (4);
            \path (4) edge (1);
            \path (4) edge (5);
            \path (5) edge (6);
            \path (6) edge (7);
            \path (7) edge (1);
        \end{tikzpicture}
        \caption{
            A graph $G$
        }
        \label{fig: just house with garage}
    \end{subfigure}
    ~
    \begin{subfigure}{0.3\textwidth}
        \centering
        \begin{tikzpicture}[scale=1.0,
            every node/.style={circle,thick,draw,inner sep=2pt},
            every edge/.style={draw,thick}]
            \node (1) at (   0,  1) {1};
            \node (2) at (   1,  1) {2};
            \node (3) at (   1,  0) {3};
            \node (4) at (   0,  0) {4};
            \node (5) at (  -1,  0) {5};
            \node (6) at (  -1,  1) {6};
            \node (7) at (-0.5,1.75) {7};
            \path (2) edge (3);
            \path (3) edge (4);
            \path (4) edge (5);
            \path (5) edge (6);
            \path (6) edge (7);
            \path (7) edge (1);
        \end{tikzpicture}
        \caption{A spanning tree $T$ of $G$}
        \label{fig: spanning tree of house with garage}
    \end{subfigure}
    \begin{subfigure}{0.35\textwidth}
        \centering
        \begin{tikzpicture}[scale=1.0,
            every node/.style={circle,thick,draw,inner sep=2pt},
            every edge/.style={draw,thick,-latex}]
            \node (1) at (   0,  1) {1};
            \node (2) at (   1,  1) {2};
            \node (3) at (   1,  0) {3};
            \node (4) at (   0,  0) {4};
            \node (5) at (  -1,  0) {5};
            \node (6) at (  -1,  1) {6};
            \node (7) at (-0.5,1.75) {7};
            \path (2) edge (3);
            \path (4) edge (3);
            \path (4) edge (5);
            \path (6) edge (5);
            \path (6) edge (7);
            \path (1) edge (7);
        \end{tikzpicture}
        \caption{
            $\digraph{T}$ with canonical orientation
        }
        \label{fig: canonical tree of house with garage}
    \end{subfigure}
    \caption{A running example}
    \label{fig: house with garage example}
\end{figure}
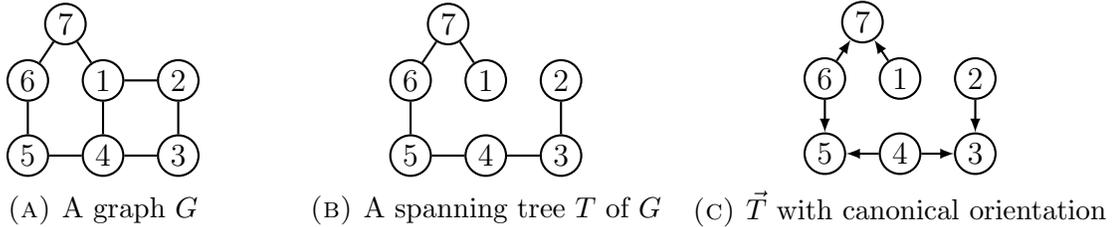

    Consider, for example, in the graph $G$ shown in \Cref{fig: just house with garage},
    which will serve as a running example.
    The subgraph $T$ in \Cref{fig: spanning tree of house with garage}
    is a spanning tree.
    In \Cref{fig: canonical tree of house with garage},
    we choose a canonical orientation and produce $\digraph{T}$.
    If we arrange its arcs into the ordered list
    \begin{equation}\label{equ: house with garage tree vector}
        ( 2 \to 3,\; 4 \to 3,\; 4 \to 5,\; 6 \to 5,\; 6 \to 7,\; 1 \to 7 ),
    \end{equation}
    then the oriented fundamental cycle
    $1 \to 2 \to \cdots \to 7 \to 1$,
    induced by $e = \{1,2\}$,
    can be expressed as the incidence vector
    $( +1, -1, +1, -1, +1, -1 )^\top$.
    
Similarly, fixing an ordering of the arcs,
a cut-set defined by a cut can be expressed as the vector with
entries in $\{ -1, 0, +1 \}$ indicating the direction in which
each arc goes across the partition (0 for not crossing the partition).
In the example shown in \Cref{fig: canonical tree of house with garage},
with respect to the cut $\nodes(G) = \{ 3,5,7 \} \cup \{ 1,2,4,6 \}$,
the cut set, which include arcs listed in \eqref{equ: house with garage tree vector}, 
can be encoded as $(-1,-1,-1,-1,-1,-1)$.
as they all go from $\{1,2,4,6\}$ to $\{3,5,7\}$.
    
    \begin{theorem}\label{thm: facet parametrization}
        For a maximal bipartite subgraph $B \le G$
        with partition $V_+ \cup V_- = \nodes(B)$,
        let $T$ be a spanning tree of $B$
        and $\digraph{T}$ be the corresponding digraph with canonical orientations
        (see \Cref{rmk: canonical orientation}).
        There is a bijection between the facets
        $\{ F \in \facets(\red{\adjp}_G) \mid G_F = B \}$
        and the set of cut-set vectors $\mathbf{d} \in \{ \pm 1 \}^n$ 
        of $T$ with respect to the cut $(V_+,V_-)$ satisfying the system
        \begin{equation}\label{equ: cycle equation}
            \begin{cases}
                \mathbf{c}_{\digraph{T}}(e)^\top \mathbf{d} &= \pm 1 
                \quad\text{for } e \in \edges(B) \setminus \edges(T) \\
                \mathbf{c}_{\digraph{T}}(e)^\top \mathbf{d} &= \phantom{\pm} 0
                \quad\text{for } e \in \edges(G) \setminus \edges(B)
            \end{cases}
        \end{equation}
    \end{theorem}
    
    Note that the $\pm$ sign is the consequence of the inherent ambiguity
    in the orientation assignment for a given fundamental cycle.
    
    \begin{proof}
        Since the maximal bipartite subgraph $B$ must span $G$,
        $T$ also spans $G$, 
        and $(V_+,V_-)$ is a partition of $\nodes(G)$.
        Recall that $\digraph{T}$ (see \Cref{rmk: canonical orientation}) has the arc set
        \[
            \edges(\digraph{T}) =
            \{ (i,j) \mid \{ i, j \} \in \edges(T) ,\, i \in V_- ,\, j \in V_+ \}.
        \]
        Let $\mathbf{d} \in \{ \pm 1 \}^n$ be a cut-set vector satisfying the system of equations
        \eqref{equ: cycle equation},
        and let $D = \operatorname{diag}(\mathbf{d})$.
        Then there exists a vector $\boldalpha \in \R^N$ such that
        $ \boldalpha^\top Q(\digraph{T}) = - \mathbf{d} $,
        i.e.
        \[ 
            \boldalpha^\top Q(\digraph{T}) D = - \boldone.
        \]
        Define $\Delta'$ to be the set of points that are the columns of 
        $Q(\digraph{T}) \operatorname{diag}(\mathbf{d})$. 
        Then $G_{\Delta'} = T$ and thus $\Delta'$ is a codimension-1 simplex in $\bar{\adjp}_G$.
        Moreover, the above equations show that
        \[
            \inner{ \boldx }{ \boldalpha } = -1
            \quad\text{for any } \boldx \in \Delta'.
        \]
        
        If $T \ne B$, then each edge $e \in \edges(B) \setminus \edges(T)$
        determines a fundamental cycle with respect to $T$ represented by
        the vector $\mathbf{c} = \mathbf{c}_T(e)$ such that $e$ corresponds to a point $\boldx = - Q(\digraph{T}) \mathbf{c}$.
        By assumption \eqref{equ: cycle equation}, $\mathbf{c}^\top \mathbf{d} = \pm 1$. 
        Therefore,
        \begin{align*}
            \boldalpha^\top \boldx 
            &= \boldalpha^\top (- Q(\digraph{T}) \mathbf{c}) 
            = - \boldalpha^\top Q(\digraph{T}) D D \mathbf{c} 
            = \boldone^\top D \mathbf{c} 
            = \mathbf{d}^\top \mathbf{c} = \pm 1.
        \end{align*}
        
        Similarly, for any $e \in \edges(G) \setminus \edges(B)$,
        the fundamental cycle with respect to $T$ is represented by a 
        vector $\mathbf{c} = \mathbf{c}_T(e)$ such that $e$ corresponds to a point
        $ \boldx = - Q(\digraph{T}) \mathbf{c}$.
        By assumption \eqref{equ: cycle equation}, $\mathbf{c}^\top \mathbf{d} = 0$. 
        Following the same calculations above,
        \begin{align*}
            \boldalpha^\top \boldx 
            &= \boldalpha^\top (- Q(\digraph{T}) \mathbf{c}) 
            = - \boldalpha^\top Q(\digraph{T}) D D \mathbf{c} 
            = \boldone^\top D \mathbf{c} 
            = \mathbf{d}^\top \mathbf{c} = 0.
        \end{align*}
        
        We have shown $\inner{ \boldx }{ \boldalpha }$ is $-1$
        for all $\boldx$ in the codimension-1 simplex $\Delta'$ of $\bar{\adjp}_G$
        and $\pm 1$ or $0$ for any non-interior point in $\bar{\adjp}_G$.
        Therefore, $\boldalpha$ defines a unique facet
        $\operatorname{aff}(\Delta') \cap \bar{\adjp}_G \in \facets(\adjp_G)$.
        That is, each solution $\mathbf{d}$ to the system \eqref{equ: cycle equation},
        determines a unique facet of $\bar{\adjp}_G$.
        
        Conversely, any facet $F \in \facets(\bar{\adjp}_G)$  such that $G_F = B \ge T$
        must contain the subset $\{ d_1 \boldx_1, \dots, d_n \boldx_n \}$
        for some $\mathbf{d} = (d_1,\dots,d_n) \in \{ \pm 1 \}^n$,
        where $\{ \boldx_1,\dots,\boldx_n \} = \bar{\adjp}_{\digraph{T}}$.
        The vector $\mathbf{d}$ is uniquely determined by the choice of $\digraph{T}$.
        By \Cref{thm: facets have half and half},
        $\mathbf{d}$ must satisfy the equations in \eqref{equ: cycle equation}.
        That is, each facet of $\bar{\adjp}_G$ corresponds to 
        a unique solution $\mathbf{d}$ to the system \eqref{equ: cycle equation}.
    \end{proof}
    
    D'Alì, Delucchi, and Michałek showed that 
    for a connected bipartite graph $G$, 
    $|\facets(\adjp_G)|$ is bounded by $2^{N-1}$
    \cite[\rev{Corollary 33}]{DAliDelucchiMichalek2022Many}.
    It is then noted that this upper bound no longer holds when the graph is not bipartite.
    From the above proof, we can derive a refinement of this result:
    this upper bound always holds for the number of facets 
    in an equivalence class of facets associated with a given facet subgraph.
    
    \begin{corollary}[\rev{A refinement of Corollary 33 of Ref.~\cite{DAliDelucchiMichalek2022Many}}]
    \label{cor: facet class bound}
        For a facet $F \in \facets(\adjp_G)$,
        \[
            | \, \{ F' \in \facets(\adjp_G) \mid G_{F'} = G_F \} \, |
            \le 2^{N-1}.
        \]
    \end{corollary}
    
    With this, we can derive an upper bound for the total number of facets,
    \rev{which is a generalization of \cite[Corollary 33]{DAliDelucchiMichalek2022Many} to all connected graphs}.
    
    \begin{corollary}\label{cor: facet total bound}
        If $\beta$ is the number of maximal bipartite subgraphs of $G$, then
        \[
            | \, \facets(\adjp_G) \, |
            \le \beta \cdot 2^{N-1}.
        \]
    \end{corollary}
    
    When $G$ is bipartite, the only facet subgraph
    (the unique maximal bipartite subgraph) is $G$ itself,
    and the result reduces to the previously established upper bound
    \cite[\rev{Corollary 33}]{DAliDelucchiMichalek2022Many}.

\subsection{Properties of faces and their face subgraphs}\label{sec: properties}

We now establish connections between the geometric properties of faces of $\adjp_G$
and the graph-theoretical properties of their corresponding face subgraphs.
Recall that the \emph{cyclomatic number} of a graph $G$ is the minimum number 
of edges that can be deleted from $G$ such that the resulting graph is acyclic.

\begin{theorem}\label{thm: face properties}
    For a proper face $F$ of $\adjp_G$,
    \begin{enumerate}[label=(\roman*)]
        \item \label{part: independent iff forest}
            $F$ is independent if and only if $G_F$ is a forest;
            
        \item \label{part: circuit iff chordless cycle}
            $F$ is a circuit if and only if $G_F$ is a chordless cycle;
            
        \item\label{part: face dimension}
            $\dim(F) = |\nodes(G_F)| - k - 1$ where $k$ is the number of connected components in $G_F$;
            
        \item\label{part: face corank} $\corank(F)$ is the cyclomatic number of $G_F$.
    \end{enumerate}
\end{theorem}

If $G_F$ is spanning, its Betti numbers are the codimension and corank of $F$.
Some properties have been studied in different context.
E.g., part \ref{part: face dimension} was established in
\cite[Lemma 1]{GordonPetrov2017}.

\begin{proof}
    \hfill
    \begin{enumerate}[label=(\roman*),wide, labelwidth=!, itemindent=0em]
    \setlength{\itemsep}{1ex}
    
    \item 
        Let $F'$ be a facet of $\adjp_G$ containing $F$.
        If $G_F$ contains a cycle with edges 
        $i_1 \leftrightarrow \cdots \leftrightarrow i_m \leftrightarrow i_1$,
        then this cycle is also contained in $G_{F'}$.
        By \Cref{thm: facets have half and half}, 
        it must be an even cycle,
        and there exist $\lambda_1,\dots,\lambda_m \in \{ \pm 1 \}$ with 
        $\sum_{j=1}^m \lambda_j = 0$ such that
        $\lambda_j (\bolde_{i_j} - \bolde_{i_{j+1}}) \in F$ for all $j$.
        This gives us the affine dependence relation
        \[
            \sum_{j=1}^m \lambda_j (\lambda_j (\bolde_{i_j} - \bolde_{i_{j+1}})) =
            \sum_{j=1}^m \rev{ (\bolde_{i_j} - \bolde_{i_{j+1}} ) } =
            \boldzero
        \]
        with the coefficients $\lambda_1,\dots,\lambda_m$.
        Therefore $F$ itself cannot be independent.
        
        Conversely, if $G_F$ is a forest, 
        then the incidence matrix $Q(\digraph{G}_F)$,
        whose columns are points in $F$, has full column rank
        \cite[Lemma 2.5]{Bapat2010Graphs} \cite{Schrijver1998Theory}.
        Therefore $F$ is independent.
            
    \item 
    
        If $F$ is a circuit, then $F$ is dependent by definition.
        By part (i), $G_F$ contains a cycle
        and the corresponding subset of points in $F$
        is dependent.
        However, the circuit $F$, being a minimal affinely dependent set,
        must be exactly this set.
        Therefore $G_F$ is exactly this cycle.
            
        Conversely, if $F$ is not a circuit,
        then either $F$ is independent or 
        $F$ contains a proper dependent subset $F'$.
        Again by part (i), 
        $G_F$ is either a forest or it contains a strictly smaller cycle,
        and thus $G_F$ is not a chordless cycle.
            
    \item 
        Let $k$ be the number of connected components in $G_{F}$.
        Then
        \[
            \dim(F) = \rank(Q(\digraph{G}_F)) - 1 = |\nodes(G_F)| - k - 1.
        \]
    
    \item 
        Let $k$ be as defined above and 
        $\mu$ be the cyclomatic number of $G_F$.
        By part (iii),
        \begin{align*}
            \mu  
            &= |\edges(G_F)| - |\nodes(G_{F})| + k 
            = |F| - (\dim(F) + 1)
            = \corank(F). \qedhere
        \end{align*}
    \end{enumerate}
\end{proof}

\begin{remark}
    \Cref{thm: face properties} 
    highlights the connection through which
    independent faces correspond to forests,
    dependent faces correspond to cyclic graphs,
    and circuit faces correspond to chordless cycles.
    The precise description emerges from the viewpoint of matroid theory
    \cite{DelucchiHoessly}.
\end{remark}

Combining \Cref{thm: facets are max bipartite,thm: face properties}
part \ref{part: independent iff forest},
we get a simple alternative proof to the fact that
$\adjp_G$ is simplicial (i.e., all of its facets are simplices)
if and only if all maximal bipartite subgraphs of $G$ are trees,
which was first established \rev{by Higashitani}
in \cite[Corollary 2.3]{Higashitani2015Smooth2}
and has important consequences in the study of facet systems of Kuramoto equations
(see \Cref{sec: applications}) and the structure of certain metric spaces \cite{GordonPetrov2017}.

\begin{corollary}[\rev{\cite[Corollary 2.3]{Higashitani2015Smooth2} }]
    \label{cor: simplicial iff no even cycle}
    $\adjp_G$ is simplicial if and only if $G$ has no even cycles.
    \qed
\end{corollary}

\section{Case study: Joining two cycles along a shared edge}\label{sec:examples}

D'Alì, Delucchi, and Michałek showed that for a graph formed by
joining two bipartite graphs along an edge,
the number of facets of the associated symmetric edge polytope
is $\frac{1}{2} f_1 f_2$
where $f_1$ and $f_2$ are the number facets in the symmetric edge polytopes
associated with the two bipartite subgraphs respectively
\cite[\rev{Proposition 37}]{DAliDelucchiMichalek2022Many}.
In the following we explore the more general situation in which
one of the subgraph is \emph{not} bipartite.

The running example $G$, from  \Cref{fig: just house with garage},
is a non-bipartite graph formed by joining a 4-cycle and a 5-cycle along a single shared edge.
As we will calculate, the facet count described above
(\cite[\rev{Proposition 37}]{DAliDelucchiMichalek2022Many}) no longer applies,
yet \Cref{thm: facet parametrization} provides a concrete recipe
for calculating the number of facets
and describing their combinatorial structures. 

\begin{figure}[ht]
    \centering
    \begin{subfigure}{0.7\textwidth}
        \begin{tikzpicture}[scale=1.0,
            every node/.style={circle,thick,draw,inner sep=2pt},
            every edge/.style={draw,thick}]
            \node (1) at (   0,  1) {1};
            \node (2) at (   1,  1) {2};
            \node (3) at (   1,  0) {3};
            \node (4) at (   0,  0) {4};
            \node (5) at (  -1,  0) {5};
            \node (6) at (  -1,  1) {6};
            \node (7) at (-0.5,1.75) {7};
            \path (1) edge (2);
            \path (2) edge (3);
            \path (4) edge (5);
            \path (5) edge (6);
            \path (6) edge (7);
            \path (7) edge (1);
        \end{tikzpicture}
        \hfill
        \begin{tikzpicture}[scale=1.0,
            every node/.style={circle,thick,draw,inner sep=2pt},
            every edge/.style={draw,thick}]
            \node (1) at (   0,  1) {1};
            \node (2) at (   1,  1) {2};
            \node (3) at (   1,  0) {3};
            \node (4) at (   0,  0) {4};
            \node (5) at (  -1,  0) {5};
            \node (6) at (  -1,  1) {6};
            \node (7) at (-0.5,1.75) {7};
            \path (1) edge (2);
            \path (3) edge (4);
            \path (4) edge (5);
            \path (5) edge (6);
            \path (6) edge (7);
            \path (7) edge (1);
        \end{tikzpicture}
        \hfill
        \begin{tikzpicture}[scale=1.0,
            every node/.style={circle,thick,draw,inner sep=2pt},
            every edge/.style={draw,thick}]
            \node (1) at (   0,  1) {1};
            \node (2) at (   1,  1) {2};
            \node (3) at (   1,  0) {3};
            \node (4) at (   0,  0) {4};
            \node (5) at (  -1,  0) {5};
            \node (6) at (  -1,  1) {6};
            \node (7) at (-0.5,1.75) {7};
            \path (2) edge (3);
            \path (3) edge (4);
            \path (4) edge (5);
            \path (5) edge (6);
            \path (6) edge (7);
            \path (7) edge (1);
        \end{tikzpicture}
        \caption{Corank-0 facet subgraphs of $G$}
        \label{fig: house with garage corank-0 graphs}
    \end{subfigure}
    \vspace{3ex}
    
    \begin{subfigure}{0.9\textwidth}
        \begin{tikzpicture}[scale=1.0,
            every node/.style={circle,thick,draw,inner sep=2pt},
            every edge/.style={draw,thick}]
            \node (1) at (   0,  1) {1};
            \node (2) at (   1,  1) {2};
            \node (3) at (   1,  0) {3};
            \node (4) at (   0,  0) {4};
            \node (5) at (  -1,  0) {5};
            \node (6) at (  -1,  1) {6};
            \node (7) at (-0.5,1.75) {7};
            \path (1) edge (2);
            \path (2) edge (3);
            \path (3) edge (4);
            \path (4) edge (1);
            \path (4) edge (5);
            \path (5) edge (6);
            \path (6) edge (7);
        \end{tikzpicture}
        \hfill
        \begin{tikzpicture}[scale=1.0,
            every node/.style={circle,thick,draw,inner sep=2pt},
            every edge/.style={draw,thick}]
            \node (1) at (   0,  1) {1};
            \node (2) at (   1,  1) {2};
            \node (3) at (   1,  0) {3};
            \node (4) at (   0,  0) {4};
            \node (5) at (  -1,  0) {5};
            \node (6) at (  -1,  1) {6};
            \node (7) at (-0.5,1.75) {7};
            \path (1) edge (2);
            \path (2) edge (3);
            \path (3) edge (4);
            \path (4) edge (1);
            \path (4) edge (5);
            \path (5) edge (6);
            \path (7) edge (1);
        \end{tikzpicture}
        \hfill
        \begin{tikzpicture}[scale=1.0,
            every node/.style={circle,thick,draw,inner sep=2pt},
            every edge/.style={draw,thick}]
            \node (1) at (   0,  1) {1};
            \node (2) at (   1,  1) {2};
            \node (3) at (   1,  0) {3};
            \node (4) at (   0,  0) {4};
            \node (5) at (  -1,  0) {5};
            \node (6) at (  -1,  1) {6};
            \node (7) at (-0.5,1.75) {7};
            \path (1) edge (2);
            \path (2) edge (3);
            \path (3) edge (4);
            \path (4) edge (1);
            \path (4) edge (5);
            \path (6) edge (7);
            \path (7) edge (1);
        \end{tikzpicture}
        \hfill
        \begin{tikzpicture}[scale=1.0,
            every node/.style={circle,thick,draw,inner sep=2pt},
            every edge/.style={draw,thick}]
            \node (1) at (   0,  1) {1};
            \node (2) at (   1,  1) {2};
            \node (3) at (   1,  0) {3};
            \node (4) at (   0,  0) {4};
            \node (5) at (  -1,  0) {5};
            \node (6) at (  -1,  1) {6};
            \node (7) at (-0.5,1.75) {7};
            \path (1) edge (2);
            \path (2) edge (3);
            \path (3) edge (4);
            \path (4) edge (1);
            \path (5) edge (6);
            \path (6) edge (7);
            \path (7) edge (1);
        \end{tikzpicture}
        \caption{Corank-1 facet subgraphs of $G$}
        \label{fig: house with garage corank-1 graphs}
    \end{subfigure}
    \caption{Facet subgraphs (maximal bipartite subgraphs) of $G$}
    \label{fig: house with garage facet graphs}
\end{figure}
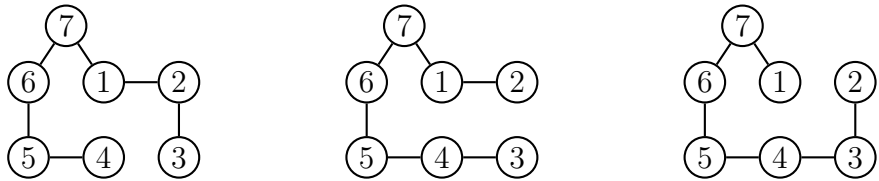
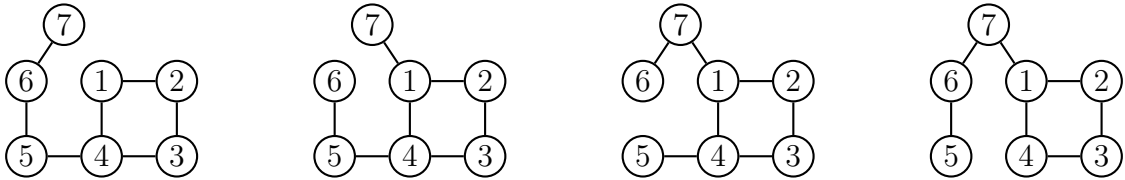

\Cref{fig: house with garage facet graphs} shows the seven maximal bipartite subgraphs of $G$.
As established in \Cref{thm: facets are max bipartite},
they correspond to the five equivalence classes of facets in $\facets(\adjp_G)$.
Among these subgraphs, three of them are trees
(\Cref{fig: house with garage corank-0 graphs}),
and, according to \Cref{thm: face properties} part \ref{part: independent iff forest},
they correspond to corank-0 (simplicial) facets.
The other four each contain a unique 4-cycle
(\Cref{fig: house with garage corank-1 graphs}),
and, according to \Cref{thm: face properties} part \ref{part: face corank},
they correspond to corank-1 facets.

As shown in \Cref{fig: house with garage example}, we pick a corank-0 facet subgraph, a spanning tree of this subgraph, and an assignment of edge orientations.
Up to a recording of the edges,
the two fundamental cycles can be expressed as vectors
$[+1,-1,+1,-1,+1,-1]$ and $[0,0,+1,-1,+1,-1]$.
Therefore the defining equation \eqref{equ: cycle equation} 
in \Cref{thm: facet parametrization}
for the parametrization of the facets in this equivalence class is
\[
    \begin{bmatrix}
       +1 & -1 & +1 & -1 & +1 & -1 \\
        0 &  0 & +1 & -1 & +1 & -1     
    \end{bmatrix}
    \mathbf{d}
    =
    \begin{bmatrix}
       0 \\ 0
    \end{bmatrix}
\]
in the unknowns $\mathbf{d} = (d_1,\dots,d_6)^\top \in \{ \pm 1 \}^6$.
This equation is equivalent to
\[
    \begin{bmatrix}
       +1 & -1 &  0 &  0 &  0 &  0 \\
        0 &  0 & +1 & -1 & +1 & -1     
    \end{bmatrix}
    \mathbf{d}
    =
    \begin{bmatrix}
       0 \\ 0
    \end{bmatrix}.
\]
From this we can see that $(d_1,d_2)$ and $(d_3,d_4,d_5)$
can be described independently,
and there are two possible choices for $(d_1,d_2)$,
namely $(+1,+1)$ and $(-1,-1)$.
Similarly, there are six possible choices for $(d_3,d_4,d_5,d_6)$:
\begin{align*}
    &(+1,+1,+1,+1) &
    &(+1,+1,-1,-1) &
    &(+1,-1,-1,+1) \\
    &(-1,-1,+1,+1) &
    &(-1,-1,-1,-1) &
    &(-1,+1,+1,-1) .
\end{align*}
Altogether, there are 12 possible choice for the vector $\mathbf{d} \in \{ \pm 1 \}^6$
for the equation \eqref{equ: cycle equation}.
These produce 12 distinct facets in the equivalence class of facets
corresponding to a corank-0 maximal bipartite subgraph of $G$
shown in \Cref{fig: house with garage corank-0 graphs}.

\begin{wrapfigure}[10]{r}{0.24\textwidth}
    \centering
    \begin{tikzpicture}[scale=1.10,
        every node/.style={circle,thick,draw,inner sep=2pt},
        every edge/.style={draw,thick,-latex}]
        \node (1) at (   0,  1) {1};
        \node (2) at (   1,  1) {2};
        \node (3) at (   1,  0) {3};
        \node (4) at (   0,  0) {4};
        \node (5) at (  -1,  0) {5};
        \node (6) at (  -1,  1) {6};
        \node (7) at (-0.5,1.75) {7};
        \path (1) edge (2);
        \path[densely dotted] (3) edge (2);
        \path (3) edge (4);
        \path (1) edge (4);
        \path (5) edge (4);
        \path (5) edge (6);
        \path (7) edge (6);
    \end{tikzpicture}
    \caption{
        A corank-1 facet subgraph
        and its spanning tree.
    }
    \label{fig: house with garage corank-1 example}
\end{wrapfigure}
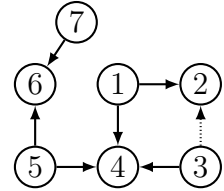
Similarly, we pick a facet subgraph $G_F$ of corank 1
in \Cref{fig: house with garage corank-1 graphs},
a spanning tree $T$ of this subgraph,
and a canonical assignment of edge orientations $\digraph{T}$
shown in \Cref{fig: house with garage corank-1 example},
where the dotted edge is in $\digraph{G}_F$ but not in $\digraph{T}$.
With respect to this choice,
and up to a reordering of the edges in $\digraph{T}$,
the two fundamental cycles can be expressed as vectors
$[-1,-1,+1,0,0,0]$ and $[0,0,+1,-1,+1,-1]$.
Therefore the fundamental cycle equations are
\[
    \begin{bmatrix}
       -1 & -1 & +1 &  0 &  0 &  0 \\
        0 &  0 & +1 & -1 & +1 & -1
    \end{bmatrix}
    \mathbf{d}
    =
    \begin{bmatrix}
       \pm 1 \\
       0
    \end{bmatrix}.
\]
Through direct calculations we can see there are 18 solutions
for $\mathbf{d} \in \{ \pm 1 \}^6$
corresponding to the 18 facets in the equivalence class.

The same argument can be applied to each of the four corank-1 facet subgraphs
in \Cref{fig: house with garage corank-1 graphs}.
Therefore there are 72 corank-1 facets in $\facets(\adjp_G)$.
All together there are $36 + 72 = 108$ facets in $\facets(\adjp_G)$.
Among them, 36 facets are simplicial
and the remaining 72 facets are (affinely) dependent and of corank 1.

The calculation shown in this concrete example can be easily generalized to
graphs formed by joining an even cycle and an odd cycle along a shared edge.
The counting argument involved makes use of the following elementary formulas from combinatorics.

\begin{lemma}\label{lem: integer solutions 1}
    For a positive integer $n$,
    there are exactly $\binom{2n}{n}$ distinct choices of vectors 
    $\mathbf{d} \in \{ \pm 1 \}^{2n}$ that satisfy the equation
    \[
        \begin{tikzpicture}[ decoration=brace, baseline=(current  bounding  box.center) ]
            \matrix (m) [
                matrix of math nodes,
                left delimiter=[,
                right delimiter={]}
                ] {
                   +1 & \cdots & +1 & -1 & \cdots & -1 \\
                };
            \draw[decorate,transform canvas={yshift=0.5em},thick] 
                (m-1-1.north) -- node[above=2pt] {$n$} (m-1-3.north);
            \draw[decorate,transform canvas={yshift=0.5em},thick] 
                (m-1-4.north) -- node[above=2pt] {$n$}     (m-1-6.north);
            \node [right = 1.5ex of m] (d) 
                {$\mathbf{d} = 0.$};
        \end{tikzpicture}
    \]
\end{lemma}

\begin{lemma}\label{lem: integer solutions 2}
    For a positive integer $n$,
    there are exactly $\binom{2n}{n-1}$ distinct choices of vectors 
    $\mathbf{d} \in \{ \pm 1 \}^{2n}$ that satisfy the equation
    \[
        \begin{tikzpicture}[ decoration=brace, baseline=(current  bounding  box.center) ]
            \matrix (m) [
                matrix of math nodes,
                left delimiter=[,
                right delimiter={]}
                ] {
                   +1 & \cdots & +1 & -1 & \cdots & -1 \\
                };
            \draw[decorate,transform canvas={yshift=0.5em},thick] 
                (m-1-1.north) -- node[above=2pt] {$n$} (m-1-3.north);
            \draw[decorate,transform canvas={yshift=0.5em},thick] 
                (m-1-4.north) -- node[above=2pt] {$n$}     (m-1-6.north);
            \node [right = 1.5ex of m] (d) 
                {$\mathbf{d} = 2.$};
        \end{tikzpicture}
    \]
\end{lemma}
\begin{proof}
    The solutions $\mathbf{d} = (d_1,\dots,d_{2n})^\top \in \{ \pm 1 \}^{2n}$ 
    correspond to the different ways of choosing 
    only $i$   entries among $(d_1,\dots,d_n)$
    and  $i+1$ entries among $(d_{n+1},\dots,d_{2n})$ to be $-1$
    for $i=0,\dots,n-1$.
    By applying the Vandermonde identity, the total number of possibilities is
    \[
        \sum_{i=0}^{n-1} \binom{n}{i} \binom{n}{i+1}
        = 
        \binom{2n}{n-1}.
        \qedhere
    \]
\end{proof}

\begin{proposition}
    Let $G$ be the graph formed by joining two cycles of size $2m_1$ and $2m_2+1$ respectively
    along a shared edge.
    The total number of facets of $\adjp_G$ is
    \[
        (2m_1-1) \binom{2m_1 - 2}{m_1 - 1} \binom{2m_2}{m_2}   + 
        (2m_2)   \binom{2m_1 - 1}{m_1}     \binom{2m_2}{m_2},
    \]
    and the two summands are the number of corank-0 (simplicial) facets
    and corank-1 facets of $\adjp_G$, respectively.
\end{proposition}

\begin{proof}
    First note that $G$ contains $N=2m_1 + 2m_2 -1$ nodes and $2m_1 + 2m_2$ edges.
    Since $G$ has a unique even cycle,
    by \Cref{thm: facets are max bipartite} and \Cref{thm: face properties},
    the coranks of facets of $\adjp_G$ are either 0 or 1.
    We shall count them separately.
    
    \textbf{(Corank-0)} We first count the corank-$0$ facets.
    By \Cref{thm: face properties} part \ref{part: face corank},
    the facet subgraph of $G$ associated with a corank-0 facet
    must be a spanning tree of $G$ that is also a maximal bipartite subgraph of $G$.
    There are exactly $2m_1 - 1$ such spanning trees of $G$, obtained by deleting the edge shared by the $2m_1$-cycle and the $(2m_2+1)$-cycle and by deleting exactly one additional edge of the $2m_1$-cycle (exemplified in \Cref{fig: house with garage corank-0 graphs}).
    Note that these corank-0 facet subgraphs are all paths of the same length, hence isomorphic to one another.
    Up to a re-indexing of the $n = 2m_1 + 2m_2 - 2$ edges, 
    the counting arguments for each of these spanning trees are identical. 
    It is therefore sufficient to calculate the number of facets
    associated with one such spanning tree (path) $T$ of $G$.
    
    There are exactly two edges in $\edges(G) \setminus \edges(T)$:
    the edge shared by the two cycles and another edge of the even cycle.
    Let $\digraph{T}$ be the corresponding digraph 
    resulting from the canonical choice of edge orientations
    as described in \Cref{rmk: canonical orientation}.
    The fundamental cycle $C$ induced by the first edge includes all edges of the odd cycle,
    and the other fundamental cycle $C'$ includes all edges of $\digraph{T}$.
    Therefore, equation~\eqref{equ: cycle equation} is of the form
    \begin{equation}\label{equ: cycle equation example}
        \begin{tikzpicture}[ decoration=brace, baseline=(current  bounding  box.center) ]
            \matrix (m) [
                matrix of math nodes,
                left delimiter=[,
                right delimiter={]}
                ] {
                    0 &  0 & \cdots &  0 &  0 & +1 & -1 & \cdots & +1 & -1 \\
                   +1 & -1 & \cdots & +1 & -1 & +1 & -1 & \cdots & +1 & -1 \\
                };
            \draw[decorate,transform canvas={yshift=0.5em},thick] 
                (m-1-1.north) -- node[above=2pt] {$2m_1 - 2$} (m-1-5.north);
            \draw[decorate,transform canvas={yshift=0.5em},thick] 
                (m-1-6.north) -- node[above=2pt] {$2m_2$}     (m-1-10.north);
            \node [right = 1.5ex of m] (d) 
                {$\mathbf{d} = \begin{bmatrix} 0 \\ 0 \end{bmatrix},$};
        \end{tikzpicture}
    \end{equation}
    and the number of facets whose facet subgraph is $T$ is exactly the
    total number of solutions $\mathbf{d} = (d_1,\dots,d_n)^\top \in \{ \pm 1 \}^n$
    to the above equation.
    Note that $(d_1,\dots, d_{2m_1 - 2}) \in \{ \pm 1 \}^{2m_1 - 2}$ 
    and $(d_{2m_1-1},\dots, d_n) \in \{ \pm 1 \}^{2m_2}$
    can be solved independently.
    
    By applying \Cref{lem: integer solutions 1}
    to these two groups of coordinates,
    we can see the number of solutions $\mathbf{d} \in \{ \pm 1 \}^n$
    to \eqref{equ: cycle equation example},
    being the product of solutions for
    $(d_1,\dots, d_{2m_1 - 2}) \in \{ \pm 1 \}^{2m_1 - 2}$ and 
    $(d_{2m_1-1},\dots, d_n) \in \{ \pm 1 \}^{2m_2}$ is
    $\binom{2m_1 - 2}{m_1 - 1} \binom{2m_2}{m_2}$.
    By \Cref{thm: facet parametrization}, 
    \[ 
        | \{ F \in \facets(\adjp_G) \mid G_F = T \} |
        =
        \binom{2m_1 - 2}{m_1 - 1} \binom{2m_2}{m_2}.
    \]
    Recall that there are $2m_1 - 1$ distinct spanning tree 
    that are maximal bipartite subgraphs of $G$,
    each having the same number of associated facets,
    therefore, the total number of corank-0 facets $\adjp_G$ has is
    \begin{equation}\label{equ: corank-0 facet count example}
        (2m_1 - 1)
        \binom{2m_1 - 2}{m_1 - 1}
        \binom{2m_2}{m_2}.
    \end{equation}
    
    \textbf{(Corank-1)}
    Now we count the corank-$1$ facets.
    Since $G$ contains a unique even cycle, 
    by \Cref{thm: facets are max bipartite} and 
    \Cref{thm: face properties} part \ref{part: face corank},
    the facet subgraph $G_F$ of a facet $F \in \facets(\adjp_G)$ of corank 1,
    being a maximal bipartite subgraph of $G$ of cyclomatic number 1,
    must contain this even cycle as well as all except one edge of the odd cycle.
    Therefore, there are exactly $2m_2$ distinct subgraphs of $G$ 
    corresponding to facets of corank 1,
    as exemplified in \Cref{fig: house with garage corank-1 graphs}.
    As in the previous case, 
    since all of these corank-1 facet subgraphs have the same fundamental cycles,
    up to a re-indexing of edges, 
    the counting arguments for each of these subgraph are identical, 
    and it is therefore sufficient to calculate the number of facets
    associated with a single corank-1 subgraph $B < G$.
    
    Fix any spanning tree $T < B$, and let $\digraph{T}$ be the corresponding digraph 
    resulting from the canonical choice of edge orientations.
    There is only one edge in $\edges(B) \setminus \edges(T)$,
    and the corresponding fundamental cycle involves all edges
    of the even cycle (as shown in \Cref{fig: house with garage corank-1 example}).
    Up to a re-indexing of the edges, the corresponding fundamental cycle vector
    can be expressed as $[+1,-1,+1,\ldots,-1,+1,0,\ldots,0]$
    with the last $2m_2-1$ coordinates being zero.
    
    Similarly, there is a unique edge in $\edges(G) \setminus \edges(B)$,
    and the corresponding fundamental cycle involves all edges
    of the $(2m_2+1)$-cycle.
    Its fundamental cycle vector can be expressed as 
    $[0,\ldots,0,+1,-1,+1,\ldots,-1,+1]$
    with the first $2m_1-2$ coordinates being zero.
    Therefore, equation~\eqref{equ: cycle equation} is of the form
    \begin{equation}\label{equ: cycle equation example corank-1}
        \begin{tikzpicture}[ decoration=brace, baseline=(current  bounding  box.center) ]
            \matrix (m) [
                matrix of math nodes,
                left delimiter=[,
                right delimiter={]}
                ] {
                   +1 & \cdots & -1 & +1 &  0 & \cdots &  0 &  0 \\
                    0 & \cdots &  0 & +1 & -1 & \cdots & +1 & -1 \\
            };
            \draw[decorate,transform canvas={yshift=0.5em},thick] 
                (m-1-1.north) -- node[above=2pt] {$2m_1 - 2$} (m-1-3.north);
            \draw[decorate,transform canvas={yshift=0.5em},thick] 
                (m-1-4.north) -- node[above=2pt] {$2m_2$}     (m-1-8.north);
            \node [right = 1.5ex of m] (d) 
                {$\mathbf{d} = \begin{bmatrix} \pm 1 \\ 0 \end{bmatrix},$};
        \end{tikzpicture}
    \end{equation}
    and its solutions $\mathbf{d} = (d_1,\dots,d_n)^\top \in \{ \pm 1 \}^n$ 
    are in one-to-one correspondence with the facets whose facet subgraph is $B$.
    
    We can solve for
    $(d_{2m_1 -1}, d_{2m_1},\dots,d_n) \in \{ \pm 1 \}^{2m_2}$,
    subjects to the constraint $d_{2m_1 -1} - d_{2m_1},\dots, +d_n = 0$,
    independently.
    Following the counting argument from the previous case,
    we can verify that there are exactly $\binom{2m_2}{m_2}$
    choices for $(d_{2m_1 -1}, d_{2m_1},\dots,d_n) \in \{ \pm 1 \}^{2m_2}$.
    
    The choices of $(d_1,\dots,d_{2m_1-2})$, however,
    depend on the value of $d_{2m_1 - 1}$,
    since they are related by the equation
    $d_1 - d_2 + \cdots - d_{2m_1-2} + d_{2m_1-1} = p$,
    where $p \in \{ \pm 1\}$.
    We consider the two cases depending on the sign of $d_{2m_1-1}/p$.
    
    If $p = d_{2m_1-1}$, the equation is equivalent to
    $d_1 - d_2 + \cdots - d_{2m_1-2} = 0$,
    and there are exactly $\binom{2m_1-2}{m_1-1}$ distinct choices for
    $(d_1,\dots,d_{2m_1-2}) \in \{ \pm 1 \}^{2m_1-2}$.
    If $p = -d_{2m_1-1}$, the equation is equivalent to
    $d_1 - d_2 + \cdots - d_{2m_1-2} = 2 (-d_{2m_1-1})$,
    and, by \Cref{lem: integer solutions 2}, 
    there are exactly $\binom{2m_1-2}{m_1-2}$ distinct choices for
    $(d_1,\dots,d_{2m_1-2}) \in \{ \pm 1 \}^{2m_1-2}$.
    In total, the number of possibilities for 
    $(d_1,\dots,d_{2m_1-2}) \in \{ \pm 1 \}^{2m_1-2}$ is
    \[
        \binom{2m_1-2}{m_1-1} + \binom{2m_1-2}{m_1-2}
        = \binom{2m_1-1}{m_1-1}
        = \binom{2m_1-1}{m_1}.
    \]
    
    Therefore, the total number of distinct choices of
    $\mathbf{d} \in \{ \pm 1 \}^n$ that satisfy~\eqref{equ: cycle equation example corank-1}
    is
    \[
        \binom{2m_1 - 1}{m_1}
        \binom{2m_2}{m_2}.
    \]
    
    This number is also the number of facets of $\adjp_G$ whose facet subgraph is $B$,
    i.e.,
    \[
        | \{ F \in \facets(\adjp_G) \mid G_F = B \} |
        =
        \binom{2m_1 - 1}{m_1}
        \binom{2m_2}{m_2}.
    \]
    Recall that there are $2m_2$ corank-1 facet subgraphs.
    Therefore, the total number of corank-1 facets $\adjp_G$ has is
    \[
        2 m_2
        \binom{2m_1 - 1}{m_1}
        \binom{2m_2}{m_2},
    \]
    which completes the proof.
\end{proof}

Note that this proof is constructive
in the sense that the facets, encoded as cut-set vectors,
can be enumerated as solutions to
\eqref{equ: corank-0 facet count example} and
\eqref{equ: cycle equation example corank-1}.
    
We conclude with an alternative formulation for the facet count provided by
the proposition above,
similar to the result established in \cite[\rev{Proposition 37}]{DAliDelucchiMichalek2022Many}.
\rev{Ohsugi and Shibata showed that}
for an even cycle $C_{2k}$,
the number of facets of $\adjp_{C_{2k}}$ is $\binom{2k}{k}$~\cite{OhsugiShibata2012Smooth}.
Using this formula, we can relate the facet count presented above
and the facet counts for symmetric edge polytopes associated with even cycles.

\begin{corollary}
    Let $G$ be the graph formed by joining two cycles of size $2m_1$
    and $2m_2 + 1$ respectively along a shared edge.
    Then
    \[
        | \facets(\adjp_G) | =
        \frac{m_1 + 2 m_2}{2}
        \,
        f_{C_{2m_1}} f_{C_{2m_2}},
    \]
    where $f_{C_{2m_1}}$ and $f_{C_{2m_2}}$ are the number of facets
    $\adjp_{C_{2m_1}}$ and $\adjp_{C_{2m_1}}$ have respectively.
\end{corollary}



\section*{Acknowledgements}

This project is motivated by a series of questions posed by
Tien-Yien Li in his 2013 lecture on solving polynomial systems.
It is also an extension of a discussion the authors had with 
Alessio D'Alì, Emanuele Delucchi, and Mateusz Michałek.
We would like to thank the anonymous referees and Thomas Zaslavsky for their many helpful comments which improved the content and exposition within this work.

\bibliographystyle{abbrv}
\bibliography{library,updated-refs}

\appendix
\section{Applications to algebraic Kuramoto equations}\label{sec: applications}

Facets of $\adjp_G$ play important roles in the study of algebraic Kuramoto equations
\cite{Chen2019Directed,Chen2019Unmixing,ChenMarecekMehtaNeimerg2019Three},
which has attracted interests from electrical engineering, biology, and chemistry.
The original Kuramoto equations model the synchronization behaviors of
a network of coupled oscillators \cite{Kuramoto1975Self},
which can be represented by a weighted graph $G$ with the
nodes $\nodes(G) = \{0,\dots,n\}$ representing the oscillators,
the edges $\edges(G)$ representing the connections among the oscillators, 
and the weights $K = \{k_{ij}\}$ representing the \emph{coupling strengths} along the edges.
Each oscillator $i$ has its own natural frequency $\omega_i$.
The dynamics of the network
can be described by the differential equations
\begin{equation}\label{equ:kuramoto-ode}
    \frac{d \theta_i}{dt} =
    \omega_i -
    \sum_{j \in \mathcal{N}_G(i)} k_{ij} \sin(\theta_{i}-\theta_{j}),
    \quad\text{for } i = 0,\dots,n,
\end{equation}
where each $\theta_i \in [0,2\pi)$ is the phase angle 
that describes the status of the $i$-th oscillator,
and  $\mathcal{N}_G(i)$ is the set of its neighbors.
\emph{Frequency synchronization} 
occurs when the competing forces reach equilibrium and 
all oscillators are tuned to the same frequency,
i.e., $\frac{d\theta_i}{dt} = c$ for a common constant $c$ for all $i$.
They are precisely the solutions to the system of equations
\begin{equation}
    \omega_i - 
    \sum_{j \in \mathcal{N}_G(i)} k_{ij} 
    \sin(\theta_{i}-\theta_{j}) = c
    \quad\text{for } i = 1,\dots,n
    \label{equ:kuramoto-sin}
\end{equation}
in the variables $\theta_1,\dots,\theta_n$.
Here $\theta_0 = 0$ is fixed as the reference phase angle.
With the substitution $x_{i} := e^{\imag \theta_{i}}$ ($x_0 = 1$),
\eqref{equ:kuramoto-sin} can be transformed into the algebraic system
\begin{equation}\label{equ:kuramoto-rat}
    F_{G,i}(x_1,\dots,x_n) = 
    \omega_i - c - 
    \sum_{j \in \mathcal{N}_G(i)} 
    \frac{k_{ij}}{2\imag}
    \left(
        \frac{x_i}{x_j} - 
        \frac{x_j}{x_i}
    \right) 
    = 0
    \quad\text{for } i = 1,\dots,n.
\end{equation}
The system $F_G = (F_{G,1},\dots,F_{G,n})$ consists of $n$ Laurent polynomial
equations in the $n$ nonzero complex variables 
$\mathbf{x} = (x_1,\dots,x_n) \in (\CC^*)^n = (\CC \setminus \{ 0 \})^n$.

Considering $F_G$ as a column vector, for any nonsingular  $n \times n$ matrix $R$, 
the systems $F_G^R = R \cdot F_{G}$ and $F_{G}$ 
have the same zero set.
For generic choices of the matrix $R$,
there is no complete cancellation of the terms, 
and thus $F_G^R$ is of the form
\begin{equation}\label{equ:kuramoto-rand}
    F_{G,k}^R = c_k^R -
    \sum_{\{i,j\} \in \edges(G)}
    \left(
        a_{ijk}^R \frac{x_i}{x_j} + 
        a_{jik}^R \frac{x_j}{x_i}
    \right)
    \quad\text{for } k = 1,\dots,n,
\end{equation}
where $c_k^R$ and $a_{ijk}^R$ are the resulting nonzero coefficients
after collection of similar terms.
This is the \emph{algebraic Kuramoto system} in its \emph{unmixed form}.
To see the connections to symmetric edge polytopes more clearly,
we shall use the vector exponent notation
\[
    (x_1,\dots,x_n)^{\left[ \begin{smallmatrix} a_1 \\ \vdots \\ a_n \end{smallmatrix} \right]}
    =
    x_1^{a_1} \cdots x_n^{a_n}.
\]
We can then write \eqref{equ:kuramoto-rand} as
\[
    F_{G,k}^R(\boldx) = 
    \sum_{\bolda \in \red{\adjp}_G}
        c(\bolda) \, \boldx^{\bolda}
    \quad\text{for } k = 1,\dots,n,
\]
where the function $c : \red{\adjp}_G \to \mathbb{C}$ captures the coefficients.
That is, $\red{\adjp}_G$ is exactly the support of the 
unmixed form of the algebraic Kuramoto system.
Facets and faces, in general, of $\red{\adjp}_G$ play particularly important roles
in the study of this algebraic system.
In the following, we highlight three of them, namely the roles
in toric deformation homotopy method (\cref{sec: homotopy}), 
root counting (\cref{sec: initial forms}), 
and construction of homogeneous coordinates (\cref{sec: cox ring}).

\subsection{Toric deformation homotopy}\label{sec: homotopy}

The toric deformation homotopy for unmixed algebraic Kuramoto equations
is a specialized  \emph{polyhedral homotopy} \cite{HuberSturmfels1995Polyhedral} 
construction for locating \emph{all} complex zeros of \eqref{equ:kuramoto-rand},
which includes \emph{all} frequency synchronization configurations.
Utilizing the topological information extracted from the underlying graph,
this homotopy construction has the potential to avoid 
the computationally expensive preprocessing steps associated with polyhedral homotopy
(e.g. mixed cell computations).
In the most basic form, it is defined by
the function $H_G : \CC^n \times \CC \to \CC^n$ with 
$H_G(\mathbf{x},t) = (H_{G,1},\dots,H_{G,n})$ given by
\begin{equation}\label{equ:homotopy}
    H_{G,k} = 
    \sum_{\bolda \in \red{\adjp}_G}
        c(\bolda) \, \boldx^{\bolda}
        t^{\omega(\bolda)}
    \quad\text{for } k = 1,\dots,n,
    \quad\text{where }
    \omega(\bolda) = 
    \begin{cases}
        0 &\text{if } \bolda = \boldzero \\
        1 &\text{otherwise}.
    \end{cases}
\end{equation}
Clearly, $H_G(\boldx,1) = F_G^R(\boldx)$.
As $t$ varies between 0 and $1$ within the interval $(0,1)$, 
$H_G(\mathbf{x},t)$ represents a smooth deformation of the system $F_G^R$,
and the corresponding complex roots also vary smoothly,
forming smooth paths reaching \emph{all} complex zeros of $F_G^R$
\cite{Chen2019Directed,HuberSturmfels1995Polyhedral}.

The starting points of the smooth paths, however,
are not well defined, as the limit points of these paths
as $t \to 0$ are not contained in $(\CC^*)^n$.
Yet, with the change of variables
\[
    x_k = y_k t^{\alpha_k}
    \quad\text{for } k=1,\dots,n,
\]
where $\red{\boldalpha} = (\alpha_1,\dots,\alpha_n)$ is a normalized
inner normal of a facet $F \in \facets(\red{\adjp}_G)$, the limit points of certain paths, in the $y$-coordinates,
are exactly the $(\CC^*)$-solutions to the subsystem
\[
    0 =
    \sum_{\bolda \in F \cup \{ \boldzero \} } c(\bolda) \, \boldx^{\bolda}
    \quad\text{for } k = 1,\dots,n
\]
defined by the facet $F$.
As the pair $(F,\red{\boldalpha})$ runs through the set of facets
and their corresponding normalized inner normals,
the solutions of the subsystems of the above form
include limit points, in the $y$-coordinate, of \emph{all} paths defined by $H_{G,k} = 0$.
Therefore, explicit descriptions to the facets of $\red{\adjp}_G$
as well as their inner normals are crucially important
in bootstrapping the toric deformation homotopy method
for solving the Kuramoto equations.

\subsection{Facet systems and the root counting problem}\label{sec: initial forms}

A closely related application is the root counting problem
for the algebraic Kuramoto equations \eqref{equ:kuramoto-rand}.
It is shown that for generic choices of the coefficients,
the total complex root count for algebraic Kuramoto equations
induced by cycles and trees is exactly the adjacency polytope bound.
When there are algebraic relations among the coefficients,
the actual root count may be strictly less.
An algebraic certificate for such decrease in root count
is provided by ``face systems''.
For a positive-dimensional face $F$ of $\red{\adjp}_G$, 
the corresponding \emph{face system} of \eqref{equ:kuramoto-rand} is given by
\begin{equation}\label{equ: face init system}
    0 =
    \sum_{\bolda \in F } c(\bolda) \, \boldx^{\bolda}
    \quad\text{for } k = 1,\dots,n.
\end{equation}
By Bernshtein's Second Theorem~\cite{Bernshtein1975Number},
if \eqref{equ: face init system} has nontrivial solutions ($\CC^*$-solutions),
then the root count for \eqref{equ:kuramoto-rand}
is strictly less than the adjacency polytope bound.
Consequently, the descriptions of faces, especially facets,
of $\red{\adjp}_G$ given in \Cref{thm: facets are max bipartite,thm: facet parametrization}
provide a foundation for studying the root count of
the algebraic Kuramoto equations~\eqref{equ:kuramoto-rand}.

\subsection{Construction of homogeneous coordinates}\label{sec: cox ring}

One issue users of homotopy continuation methods have to face
is the existence of divergent paths,
e.g., solution paths defined by \eqref{equ:homotopy} 
that do not have limit points, as $t \to 1$, in the work space $\CC^n$.
A common solution is to compactify the work space.
Homogeneous coordinates, a special case of Cox's homogeneous ring
\cite{Cox1995Homogeneous},
provide the most general construction for achieving this goal.
Even though the use of homogeneous coordinates tends to introduce
a large number of auxiliary variables and hence is impractical in actual calculations,
they remain important tools for theoretical analysis of algebraic Kuramoto equations.
The foundation of homogeneous coordinates construction
is the full description of facets $\facets(\adjp_G)$ 
and their inner normals.

Let $m = | \facets(\red{\adjp}_G) |$,
and let $\boldalpha_1,\dots,\boldalpha_m$ be the primitive inner normals
for the facets of $\red{\adjp}_G$.
Define the matrix $V$ to be the $m \times n$ matrix whose rows are
$\boldalpha_1^\top,\dots,\boldalpha_m^\top$ and let $\mathbf{h} = (h_1,\dots,h_m)^\top$ be the column vector with entries
\[
    h_i = \min_{\boldx \in \red{\adjp}_G} \inner{ \boldx }{ \boldalpha_i }.
\]
The homogenization of $F^R_{G,k}$ \eqref{equ:kuramoto-rand} 
is the polynomial 
$\hat{F}^R_{G,k} (\boldy)$ in the complex variables 
$\boldy = (y_1,\dots,y_m)$ given by
\[
    \hat{F}^R_{G,k} (\boldy) = 
    \sum_{\bolda \in \red{\adjp}_G} c(\bolda) \, \boldy^{V \bolda - \mathbf{h}}.
\]
The system $\hat{F}^R_G = (\hat{F}^R_{G,1},\dots,\hat{F}^R_{G,n})$ 
represents a lifting of the system \eqref{equ:kuramoto-rand}
to a compact topological space (a compact toric variety).
Consequently, each solution path defined by the homotopy~\eqref{equ:homotopy}
has a limit point corresponding to an equivalence class of zero of $\hat{F}^R_G$
in this compact space, even if it has no limit point in $\CC^n$.

\end{document}